\documentclass{amsart}

\usepackage[margin=1.4in]{geometry}

\usepackage[shortlabels]{enumitem}

\usepackage{amssymb,amsmath,amsfonts,latexsym,tikz}
\newtheorem{theorem}{Theorem}
\newtheorem{definition}[theorem]{Definition}
\newtheorem{notation}[theorem]{Notation}
\newtheorem{proposition}[theorem]{Proposition}
\newtheorem{lemma}[theorem]{Lemma}
\newtheorem*{lemma*}{Lemma}
\newtheorem{corollary}[theorem]{Corollary}

\theoremstyle{definition}
\newtheorem{example}[theorem]{Example}
\newtheorem{remark}[theorem]{Remark}

\numberwithin{theorem}{section}

\DeclareMathOperator{\Aut}{Aut}

\def\A{\ensuremath{\mathbb{A}}}
\def\C{\ensuremath{\mathbb{C}}}
\def\F{\ensuremath{\mathbb{F}}}
\def\G{\ensuremath{\mathbb{G}}}
\def\cH{\ensuremath{\mathcal{H}}}
\def\L{\ensuremath{\mathbb{L}}}

\def\P{\ensuremath{\mathbb{P}}}
\def\Q{\ensuremath{\mathbb{Q}}} 

\def\SS{\ensuremath{\mathbb{S}}}
\def\Z{\ensuremath{\mathbb{Z}}}
\def\cC{\ensuremath{\mathcal{C}}}
\def\cM{\ensuremath{\mathcal{M}}}

\def\ocM{\ensuremath{\overline{\mathcal{M}}}}
\def\fX{\ensuremath{\mathfrak{X}}}
\def\V{\ensuremath{\mathbb{V}}}

\def\p{\ensuremath{\mathbf{p}}}

\def\Qcone{\ensuremath{Q^\mathrm{con}}}
\def\Qsplit{\ensuremath{Q^\mathrm{spl}}}
\def\Qns{\ensuremath{Q^\mathrm{nsp}}}

\def\Cns{\ensuremath{\cC^\mathrm{nsp}}}
\def\Cone{\ensuremath{\cC^\mathrm{con}}}
\def\Csplit{\ensuremath{\cC^\mathrm{spl}}}

\def\Ncone{\ensuremath{N^\mathrm{con}}}
\def\Nsplit{\ensuremath{N^\mathrm{spl}}}
\def\Nns{\ensuremath{N^\mathrm{nsp}}}

\def\red{\ensuremath{\mathrm{red}}}
\def\sing{\ensuremath{\mathrm{sing}}}

\def\Tr{\ensuremath{\mathrm{Tr}}}

\def\deg{\mathrm{deg}}
\def\Gal{\ensuremath{\mathrm{Gal}}}

\def\PGL{\ensuremath{\mathrm{PGL}}}

\def\Spec{\ensuremath{\mathrm{Spec \,}}}

\newcommand\Tstrut{\rule{0pt}{2.2ex}}         

\title[Odd cohomology vanishing on moduli spaces of stable curves]{Polynomial point counts and odd cohomology vanishing on moduli spaces of stable curves}

\author{Jonas Bergstr\"om}
\email{jonasb@math.su.se}
\author{Carel Faber}
\email{C.F.Faber@uu.nl}
\author{Sam Payne}
\email{sampayne@utexas.edu}
\dedicatory{To the memory of Bas Edixhoven}
\thanks{SP has been supported by NSF grants DMS--2001502 and DMS--2053261. Portions of this work were carried out while he was in residence at ICERM for the special thematic semester program on Combinatorial Algebraic Geometry, with support from NSF grant DMS--1439786 and the Simons Foundation Institute Grant Award ID 507536. Other portions of this work were carried out while the authors were participating in the thematic semester Moduli and Algebraic Cycles at Institut Mittag-Leffler}

\begin{document}
\renewcommand*{\thefootnote}{\arabic{footnote}}

\begin{abstract}
We compute the number of $\F_q$-points on $\ocM_{4,n}$ for $n \leq 3$ and show that it is a polynomial in $q$, using a sieve based on Hasse--Weil zeta functions.  As an application, we prove that the rational singular cohomology group $H^k (\ocM_{g,n})$ vanishes for all odd $k \leq 9$.  Both results confirm predictions of the Langlands program, via the conjectural correspondence with polarized algebraic cuspidal automorphic representations of conductor 1, which are classified in low weight. Our vanishing result for odd cohomology resolves a problem posed by Arbarello and Cornalba in the 1990s.
\end{abstract}

\maketitle

\section{Introduction} \label{sec:introduction}

In the 1990s, Arbarello and Cornalba introduced an inductive method for computing the rational singular cohomology groups of moduli spaces of stable curves, using the combinatorial stratification of the boundary and basic properties of Deligne's weight filtration.  As one  application, they showed that the rational singular cohomology groups $H^k(\ocM_{g,n})$ vanish for odd $k \leq 5$, for all $g$ and $n$ \cite{ArbarelloCornalba98}. Unlike the odd cohomology groups of Grassmannians and smooth projective toric varieties, the odd cohomology groups of moduli spaces of stable curves do not vanish in all degrees. For example $H^{11}(\ocM_{1,11}) \cong \Q^2$, and $\ocM_{1,n}$ has odd cohomology in a range of degrees greater than or equal to $11$, for $n > 11$ \cite{Getzler98}. It is also known that $H^{33}(\ocM_{g,n})$ is nonzero when $g$ is sufficiently large, for arbitrary $n$, as is $H^{13}(\ocM_{g,n})$ when $g$ is sufficiently large and $n \geq 10$ \cite[Corollary~4.7]{Pikaart95}.

\begin{theorem} \label{thm:oddvanishing}
The rational singular cohomology groups $H^k(\ocM_{g,n})$ vanish for all odd $k \leq 9$.
\end{theorem}

\noindent This affirmatively resolves the natural open problem highlighted by Arbarello and Cornalba  \cite[p.~103]{ArbarelloCornalba98}.  The examples in genus 1 show that the bound is best possible.  Here and throughout, all singular cohomology is taken with rational coefficients. 

\begin{remark}
This vanishing of odd cohomology in degrees less than 11 was expected based on conjectures from the Langlands program.  Simple motives appearing in the cohomology of smooth and proper Deligne--Mumford (DM) stacks over $\Z$, such as $\ocM_{g,n}$, up to Tate twists, should correspond to irreducible polarized algebraic cuspidal automorphic representations of $\PGL_n$ of conductor 1. See, e.g., the motivational discussion in \cite[\S1.2]{ChenevierRenard15}. However, there are no such representations in any odd weight less than 11 \cite[\S III, Remark 1]{Mestre86}.  For the full classification up to motivic weight 22, see  \cite[Theorem~F]{ChenevierLannes19}.
\end{remark}

Our proof of Theorem~\ref{thm:oddvanishing} follows the original inductive method of Arbarello and Cornalba.  In order to run the induction for $k = 9$, we must establish several new base cases.  In particular, we need to show that $H^9 (\ocM_{4,n})$ vanishes for $n \leq 3$.  We prove this using the Behrend--Grothendieck--Lefschetz trace formula for algebraic stacks 
\cite{BehrendThesis, Behrend93}
and point counting over finite fields.

\medskip

Recall that if $\fX$ is an algebraic stack and $\F_q$ is the finite field with $q$ elements, then $\fX(\F_q)$ is a finite groupoid.  By definition $$\#\fX(\F_q) := \sum_{x \in [\fX(\F_q)]} \frac{1}{\# \Aut(x)},$$ 
where $[\fX(\F_q)]$ denotes the set of isomorphism classes in the groupoid $\fX(\F_q)$.

\begin{proposition}[{\cite{BehrendThesis, Behrend93}}] \label{prop:pointcount}
Let $\fX$ be an algebraic stack.
\begin{enumerate}[(i)]
\item If $\fX$ is stratified by locally closed substacks $\fX = \fX_0 \sqcup \cdots \sqcup \fX_n$ then $$\#\fX(\F_q) = \#\fX_0(\F_q) + \cdots + \#\fX_n(\F_q).$$
\item If $\fX$ is a quotient stack $[X/G]$, where $X$ is a scheme and $G$ is a connected linear algebraic group, then $$\#\fX(\F_q) = \#X(\F_q) / \#G(\F_q).$$
\item If $\fX$ is a DM stack with coarse moduli space $X$, then $\#\fX(\F_q) = \#X(\F_q)$.
\item \label{it:trace} If $\fX$ is a DM stack over $\F_q$ and $\ell$ is prime to $q$ then
\begin{equation} \label{eq:trace}
\# \fX(\F_q) = \Tr(F_q \mid H^\bullet_{c, \text{\'et}} (\fX_{\overline \F_q}, \Q_\ell)),
\end{equation}
where $F_q$ is the geometric Frobenius endomorphism and $\Tr$ denotes the graded trace.
\end{enumerate}
\end{proposition}

\noindent Proposition~\ref{prop:pointcount} reduces the problem of point counting on many stacks, including all of those that we consider in this paper, to point counting on algebraic varieties.

We will be especially interested in stacks with the  \emph{polynomial point count} property, i.e., the property that $\#\fX(\F_q)$ is a polynomial in $q$. Let $\fX$ be a DM stack that is smooth, proper and of pure relative dimension over $\Z$, such as $\ocM_{g,n}$. By \cite[Theorem 2.1]{vandenBogaartEdixhoven05}, the following are equivalent:
\begin{enumerate}[(i)]
\item There is a polynomial $P(t)$ with rational coefficients and a subset $S$ of primes of Dirichlet density $1$ such that $\#\fX(\F_{p^m}) = P(p^m) + o(p^{m \cdot \dim \fX/2})$, for $p \in S$ and $m\ge1$.
\item There is a polynomial $P(t)$ with positive integer coefficients such that $\#\fX(\F_q) = P(q)$ for all finite fields $\F_q$.
\item For all primes $\ell$, the $\ell$-adic \'etale cohomology $H^\bullet_{\text{\'et}}(\fX_{\overline \Q}, \Q_\ell)$ vanishes in odd degrees, and in degree $2i$ it is isomorphic to a direct sum of copies of $\Q_{\ell}(-i)$.
\end{enumerate}
When these conditions hold, the coefficient of $t^i$ is equal to the dimension of $H^{2i}_{\text{\'et}}(\fX_{\overline \Q}, \Q_\ell)$. By the \'etale-to-singular comparison theorem and the universal coefficient theorem, if $\fX$ has polynomial point counts, then the rational singular cohomology $H^\bullet(\fX_\C)$ is supported in even degrees, and $h^{2i}$ is the coefficient of $t^i$.  Assuming furthermore that the coarse moduli space $X_\Q$ is the quotient of a smooth projective variety by a finite group, as is the case for $\ocM_{g,n}$ \cite{BoggiPikaart00}, it follows also that the Hodge structure on $H^\bullet(\fX_\C)$ is pure Tate.  

\begin{theorem} \label{thm:polynomial}
For any $n \leq 3$, $\#\ocM_{4,n}(\F_q)$ is a polynomial in $q$.  More precisely: 
\begin{align*}
\#\ocM_{4}(\F_q) & =  q^9+4q^8+13q^7+32q^6+50q^5+50q^4+32q^3+13q^2+4q+1; \\
\#\ocM_{4,1}(\F_q) & =  q^{10}+6q^9+30q^8+93q^7+191q^6+240q^5+191q^4+93q^3+30q^2+6q+1;\\
\#\ocM_{4,2}(\F_q) & = q^{11}+11q^{10}+76q^9+319q^8+838q^7+1362q^6+1362q^5+
\cdots +11q+1;\\
\#\ocM_{4,3}(\F_q) & = q^{12}+21q^{11}+207q^{10}+1168q^9+3977q^8+8296q^7+10605q^6+
\cdots +21q+1.
\end{align*}
\end{theorem}

\noindent As a consequence, the cohomology of $\ocM_{4,n}$ is pure Tate for $n \leq 3$, with Poincar\'e polynomials 
\[
h_{\ocM_{4}}(t) = t^{18}+4t^{16}+13t^{14}+32t^{12}+50t^{10}+50t^8+32t^6+13t^4+4t^2+1,
\]
and so on. The elided terms for $n = 2$ and $3$ are determined by Poincar\'e duality. The polynomial $h_{\ocM_{4}}(t)$ was previously determined in \cite{BergstromTommasi07}; for $n  \in \{1, 2, 3\}$, these computations are new.

We also determine the $\SS_{n}$-equivariant (polynomial) point counts of $\ocM_{4,n}$ for $n \leq 3$; by which we mean a point count on forms of $\ocM_{4,n}$ twisted by elements of $\SS_{n}$.
Furthermore, from Theorem~\ref{thm:polynomial} and previously known results in lower genus, we determine the point counts on the open moduli space $\cM_{4,n}$ for $n \leq 3$.

\begin{theorem} \label{thm:openpolynomial}
For any $n \leq 3$, $\#\cM_{4,n}(\F_q)$ is a polynomial in $q$.  More precisely: 
\begin{align*}
\#\cM_{4}(\F_q) & =  q^9+q^8 + q^7 - q^6; \\
\#\cM_{4,1}(\F_q) & =  q^{10}+2q^9+2q^8 - q^7 - q^6 - q^2;\\
\#\cM_{4,2}(\F_q) & = q^{11}+3q^{10}+4q^9-2q^8 - 4q^7 -2q^3 -2q^2;\\
\#\cM_{4,3}(\F_q) & = q^{12}+4q^{11}+7q^{10}-4q^9-13q^8+4q^7-q^6-11q^3+2q^2 + 2q-1.
\end{align*}
\end{theorem}

\noindent These point counts are determined from Theorem~\ref{thm:polynomial} by subtracting the point counts of the boundary $\partial \cM_{g,n} := \ocM_{g,n}  \smallsetminus \cM_{g,n}$. The boundary point counts are determined (by computer-aided calculation) using the combinatorial structure of the boundary, as encoded in modular operads \cite{GetzlerKapranov98}, plus previously known ($\SS_{\bullet}$-equivariant) polynomial point count computations for smaller moduli spaces.

\begin{remark} \label{rem:eulergalois}
If we replace $q$ in Theorem~\ref{thm:openpolynomial} by the cyclotomic character $\Q_{\ell}(-1)$, then we get an equality for the $\ell$-adic compactly supported Euler characteristic of $(\cM_{4,n})_{\overline{\Q}}$ with values in the Grothendieck group of  $\ell$-adic absolute Galois representations; see \cite[Theorem 3.2]{Bergstrom08}. Explicitly, the first equality then reads, 
\[
e_c(\cM_{4} \otimes \overline{\Q},\Q_{\ell}) =  \Q_{\ell}(-9)+\Q_{\ell}(-8)+\Q_{\ell}(-7)-\Q_{\ell}(-6). 
\]
Similarly, we can translate this into an equality for the Euler characteristic $e_c((\cM_{4})_\C, \Q)$ with values in the Grothendieck group of rational mixed Hodge structures with $q$ replaced by $\Q(-1)$, the rational Tate Hodge structure of weight $2$, see \cite[Theorem 1.4]{Tommasi05} and \cite[Theorem~8]{BergstromTommasi07}. 
\end{remark}

\begin{remark}
The fact that $\ocM_{4,n}$ has polynomial point count for $n \leq 3$ is again predicted by the Langlands program via low weight classification results.  Indeed, let $\fX$ be a DM stack that is smooth and proper over $\Z$ of relative dimension less than or equal to 12.  Assuming the conjectured correspondence between motives of conductor 1 and automorphic representations, the only motives that can appear in the cohomology of $\fX$ are powers of the Tate motive $\L$, the motive $S[12]$, of weight 11, and its twist $\L S[12]$ of weight 13 (appearing exactly when $\dim_{\Z} \fX=12$ and $S[12]$ appears).
The motive $S[12]$ corresponds to the cuspidal representation $\Delta_{11}$ and contributes to the Hodge number $h^{11,0}$. Therefore, it should appear in the motive of $\fX$ if and only if $\fX_\C$ possesses a global holomorphic $11$-form.  If $\fX_\C$ is unirational, then it has no such holomorphic form and its motive should be a polynomial in $\L$.  In particular, given that $(\ocM_{4,n})_\C$ is unirational for $n \leq 15$ \cite[Theorem~7.1]{Logan03}, the motive of $\ocM_{4,n}$ should be a polynomial in $\L$ for $n \leq 3$, and hence $\ocM_{4,n}$ should have polynomial point count. Theorem~\ref{thm:polynomial} confirms this prediction unconditionally. 
\end{remark}

\begin{remark} \label{rem:computations}
Several of the polynomials in $q$ occurring in this paper have been checked against the data
provided by computer calculations for small~$q$. For~$q=2$, the starting point is the census
of genus~$4$ curves over $\F_2$ by Xarles \cite{Xarles20}; in addition, one needs to determine
the order of the automorphism group over $\F_2$ of (a representative of) each isomorphism class.
For $q=3$, we have carried out a similar census of trigonal curves of genus~$4$ over $\F_3$, with the
orders of the automorphism groups and their numbers of points over $\F_3$, $\F_9$, and $\F_{27}$.
Since the equivariant count of hyperelliptic curves is known (cf.~\S\ref{sec:equivarianthypers}),
this is enough to verify
Theorems~\ref{thm:openpolynomial}, 
\ref{thm:equiv-point-counts}, 	
and \ref{thm:M4loc} 		
for $q=2$ and $q=3$. 
In addition, Propositions
\ref{prop:Cone0},                
\ref{prop:ConeSieve},            
\ref{prop:ns0},                  
\ref{prop:NonsplitSieve},        
\ref{prop:split0},               
\ref{prop:SplitSieve},           
\ref{prop:ConeEquiv},            
\ref{prop:NonsplitEquiv},        
and \ref{prop:SplitEquiv}        
have been verified for $q=2$, by computer calculations of representatives of all trigonal curves
over $\F_2$, together with the orders of their automorphism groups as well as
their numbers of points over $\F_2$, $\F_4$, and $\F_8$, and how many of these points are singular.
Note that we have exact formulas in all of these cases, but have stated them in some cases
only up to $o(q^6)$.
Further verifications of the counts of curves on the split quadric
are provided by van Rooij \cite{vanRooij}. 

Finally, we have compared our results with computations in the tautological subring of $H^\bullet(\ocM_{4,n})$, as follows.  The program admcycles \cite{admcycles} can compute the dimension of the subspace of $H^{2i}(\ocM_{4,n})$ generated by tautological classes, assuming that the Pixton relations span all relations among the tautological generators.  The range of degrees in which such computations can be carried through is limited by computing power; we have run the program for $n = 1$, $i \leq 4$; $n = 2$, $i \leq 3$, and $n = 3, i \leq 3$.  In all of these cases, the dimensions predicted by admcycles agree with the Betti numbers given by Theorem~\ref{thm:polynomial}.  Canning and Larson have subsequently shown that $H^\bullet(\ocM_{4,n})$ is generated by tautological classes for $n \leq 6$ \cite[Theorem~1.4]{CanningLarson22}.  It follows, in particular, that Pixton's relations are complete in the range of cases noted above.


\end{remark}

\section{Inductive argument}

Here we recall the inductive method of Arbarello and Cornalba \cite{ArbarelloCornalba98}, give the proof of Theorem~\ref{thm:oddvanishing} for $k = 7$, and explain how the case $k = 9$ follows from Theorem~\ref{thm:polynomial}. 

To start, recall that excision of the boundary $\partial \cM_{g,n} := \ocM_{g,n} \smallsetminus \cM_{g,n}$ gives a long exact sequence of rational (compactly supported) cohomology groups:
\[
\cdots \to H^k_c(\cM_{g,n}) \to H^k(\ocM_{g,n}) \to H^k(\partial \cM_{g,n}) \to H^{k+1}_c(\cM_{g,n}) \to \cdots \ .
\]
Thus, whenever $H^k_c(\cM_{g,n}) = 0$, restriction to the boundary gives an injection from $H^k(\ocM_{g,n})$ to $H^k(\partial \cM_{g,n})$.  Using ``a bit of Hodge theory," i.e., basic properties of Deligne's weight filtration, Arbarello and Cornalba improve this statement as follows. 
\begin{lemma*}[{\cite[Lemma~2.6]{ArbarelloCornalba98}}]
Suppose $H^k_c(\cM_{g,n}) = 0$.  Then pullback to the normalization of the boundary gives an injection
$
H^k(\ocM_{g,n}) \hookrightarrow H^k(\widetilde{\partial \cM_{g,n}}).
$
\end{lemma*}
\noindent Since each component of the normalization of the boundary is the quotient of a product of smaller moduli spaces $\cM_{g',n'}$, as in Notation~\ref{not:order}, by a finite group, by applying the K\"unneth formula and taking invariants, one gets an inductive method for proving odd cohomology vanishing. The key is to understand vanishing of odd (compactly supported) cohomology of the open moduli spaces $\cM_{g,n}$. 

Arbarello and Cornalba obtained vanishing in a range of degrees using Poincar\'e duality and the virtual cohomological dimension (vcd) of $\cM_{g,n}$.  We will need the following improvement. 
\begin{proposition} \label{prop:improvedbound}
Assume $g \geq 1$, then:
\begin{equation} \label{eq:improvedbound}
H^k_c(\cM_{g,n}) = 0  \quad \mbox{ for } \quad \left\{ \begin{array}{ll} k < 2g & \mbox{ and } n = 0, 1; \\ k < 2g - 2 + n &\mbox{ and } n \geq 2. \end{array} \right.
\end{equation}
\end{proposition}
\begin{proof}
The bound for $n > 1$ is \cite[(2.4)]{ArbarelloCornalba98}.  For $n = 0, 1$, the proposition improves their bound by 1. The proof is essentially identical, using Poincar\'e duality and the vanishing of $H^k(\cM_{g})$ and $H^k(\cM_{g,1})$ for large $k$.  Arbarello and Cornalba used that the vcds of $\cM_{g}$ and $\cM_{g,1}$ are $4g-5$ and $4g-3$, respectively. The improvement comes from the fact that $H^{4g-5}(\cM_g)$ and $H^{4g-3}(\cM_{g,1})$ both vanish \cite{ChurchFarbPutman12, MoritaSakasaiSuzuki13}.  
\end{proof}

\begin{remark}
Since this paper was written, Wong has improved Proposition~\ref{eq:improvedbound} for $n = 2$, showing that $H^{2g}_c(\cM_{g,2}) = 0$ for $g \geq 1$ \cite[Theorem~4.1]{Wong22}.
\end{remark}

\begin{remark}
Note that the cohomology of $\cM_{g,n}$ does not always vanish in degree equal to the vcd, for $n > 1$.  Already for $g = 2$, explicit computations show that the top weight part of $H^{n+2}_c(\cM_{2,n})$ does not vanish for $4 \leq n \leq 17$ \cite{Chan21, BibbyChanGadishYun21}.  To the best of our understanding, it is not known whether $W_k H^k(\cM_{g,n})$ vanishes in degree equal to $\mathrm{vcd}(\cM_{g,n})$. Such a weaker vanishing statement would suffice for the purposes of the Arbarello-Cornalba inductive arguments.
\end{remark}

\begin{lemma*}[{\cite[Lemma~2.9]{ArbarelloCornalba98}}]
Let $k$ be an odd integer.
Suppose $H^q(\ocM_{g,n}) = 0$ for all odd $q \leq k$ and all $g$ and $n$ such that $H^q_c(\cM_{g,n}) \neq 0$. Then $H^q(\ocM_{g,n}) = 0$ for all odd $q \leq k$ and all $g$ and $n$.
\end{lemma*}

\noindent Note that, by \eqref{eq:improvedbound}, there are only finitely many $(g,n)$ such that $H^q_c(\cM_{g,n}) \neq 0$, and they are explicitly bounded.  These are the base cases required to run the induction.

\begin{proof}[Proof of Theorem~\ref{thm:oddvanishing} for $k = 7$]
We already know that $H^q(\ocM_{g,n})$ vanishes for $q \in \{1, 3, 5\}$, for all $g$ and $n$ \cite[Theorem~2.1]{ArbarelloCornalba98}.  Also, $H^\bullet(\ocM_{0,n})$ is supported in even degrees \cite{Keel92}.  It remains to check that $H^7(\ocM_{g,n})$ vanishes in the finitely many cases where $g \geq 1$ and $\max\{ 2g, 2g -2 + n \} \leq 7$.  In fact, in each of these cases, $\ocM_{g,n}$ has polynomial point count and hence has cohomology supported in even degrees. Moreover, this holds in a slightly wider range of cases that includes all cases with $g \leq 3$ and $2g - 2 + n \leq 9$.  Indeed, we know that $\ocM_{g,n}$ has polynomial point count for $g =1$ and $n \leq 10$ \cite{Getzler98}, for $g = 2$ and $n \leq 7$ \cite[\S11.2]{Bergstrom09}, and for $g = 3$ and $n \leq 5$ \cite{Bergstrom08}.  
\end{proof}

\begin{proof}[Proof that Theorem~\ref{thm:oddvanishing} for $k = 9$ follows from Theorem~\ref{thm:polynomial}]
In the proof for $k = 7$, we have seen that the cohomology of $\ocM_{g,n}$ is supported in even degrees for a range that includes all cases with $g \leq 3$ and $2g -2 + n \leq 9$.  In order to run the induction and prove Theorem~\ref{thm:oddvanishing} for $k = 9$, it remains to check that $H^9(\ocM_{4,n}) = 0$ for $n \leq 3$. Theorem~\ref{thm:polynomial} says that $\ocM_{4,n}$ has polynomial point count for $n \leq 3$, and the required vanishing follows.
\end{proof}

\begin{remark}
The improved vanishing bound in \eqref{eq:improvedbound}, using the vanishing of singular cohomology of $\cM_{g,n}$ in degree equal to the vcd for $n = 0,1$ \cite{ChurchFarbPutman12, MoritaSakasaiSuzuki13}, is essential to our proof of Theorem~\ref{thm:oddvanishing} for $k = 9$. Without this, additional arguments would be required to prove that $H^9(\ocM_5)$ and $H^9(\ocM_{5,1})$ both vanish.   The improved bound is also used in the proof for $k = 7$, as presented above.  Without it, we would also need the vanishing of $H^7(\ocM_4)$ and $H^7(\ocM_{4,1})$ to run the induction; the proof for $k = 7$ would then also depend on Theorem~\ref{thm:polynomial}.
\end{remark}

\begin{remark}
In the proof of Theorem~\ref{thm:oddvanishing} for $k = 9$, we do not use the full strength of Theorem~\ref{thm:polynomial}; we only use the vanishing of $H^9(\ocM_{4,n})$, for $n \leq 3$.  By 
Corollary~\ref{cor:approximate}
it is enough to prove that $\#\ocM_{4,n}(\F_q)$ is a polynomial in $q$ plus $o(q^{\frac{9}{2}+n})$. In particular, some of the more subtle point counting computations in \S\S\ref{sec:cone}--\ref{sec:split} are not necessary for the proof of Theorem~\ref{thm:oddvanishing}. Nevertheless, knowing the full cohomology of $\ocM_{4,n}$ with its $\SS_n$-action will be useful for future work, e.g., for the computation of the cohomology of $\ocM_{4,4}$, $\ocM_{5}$, and $\ocM_{5,1}$. Also, as explained in the introduction, the fact that $\ocM_{4,n}$ has polynomial point count for $n \leq 3$ is predicted by the Langlands program, and the additional computations that we present are needed to prove Theorem~\ref{thm:polynomial} and confirm this prediction unconditionally.  Note that the argument involved in this prediction does not extend to $n = 4$ because, a priori, a Tate twist of the motive $S[12]$ could appear in the middle degree cohomology group $H^{13}(\ocM_{4,4})$.  However, Canning and Larson have now shown that this does not occur \cite[Theorem~1.4]{CanningLarson22}
\end{remark}

\section{Approximately polynomial point counts}

Here we remark that Theorem~\ref{thm:polynomial} gives far more precise point counting information than required for the proof of Theorem~\ref{thm:oddvanishing}.  Indeed, one can show that $H^k(\ocM_{4,n})$ vanishes for $k \leq 9$ by fixing a single prime $p$, giving approximate point counts over $\F_q$ for $q = p^m$, up to $o(q^{\frac{9}{2} + n})$, and applying the following general fact. 

Fix a prime $p$ and let $\Z_p$ denote the $p$-adic integers. 

\begin{proposition} \label{prop:approximate} Let $\fX$ be a smooth and proper DM stack of relative dimension $d$ over~$\Z_p$, and let $\widetilde F_p$ be a lift of the geometric Frobenius endomorphism.  Fix an integer $s \geq d$ and assume there is a polynomial $P(t)=\sum_i P_i \, t^i$ with rational coefficients such that $\#\fX(\F_{p^m}) = P(p^m) + o(p^{ms/2})$ as $m \to \infty$. 
Then 
\begin{itemize}
\item $\dim H^i(\fX_{\C})=0$ for all odd $i  \geq s$,  
\item $\dim H^{2i}(\fX_{\C})=P_i$ for all $s/2 \leq i \leq d$, and 
\item all eigenvalues of $\widetilde F_p$ acting on $H_{\text{\'et}}^{2i}(\fX_{\overline{\Q}_p},\Q_{\ell})$, for $\ell \neq p$, are $p^i$, for $s/2 \leq i \leq d$. 
\end{itemize}
\end{proposition}
\begin{proof} Let $d_i=\dim H_{\text{\'et}}^{i}(\fX_{\overline{\Q}_p},\Q_{\ell})$, and let $\{ \alpha_{i,j}, 1 \leq j \leq d_i\} $ be the eigenvalues of $\widetilde F_p$ acting on $H_{\text{\'et}}^{i}(\fX_{\overline{\Q}_p},\Q_{\ell})$.  By comparing $H_{\text{\'et}}^{i}(\fX_{\overline{\Q}_p},\Q_{\ell})$ with $H_{\text{\'et}}^{i}(\fX_{\overline{\F}_p},\Q_{\ell})$ (\cite[Proposition 3.1]{vandenBogaartEdixhoven05}) and applying the Behrend--Grothendieck--Lefschetz trace formula \eqref{eq:trace}, we have 
\[
\sum_{i=0}^{2d} (-1)^i \sum_{1 \leq j \leq d_i} \alpha_{i,j}^m= P(p^m) + o(p^{md/2})
\] 
as $m \to \infty$. Then \cite[Lemma 4.1]{vandenBogaartEdixhoven05} tells us that $d_i = 0$ for all odd $i \geq s$, and all eigenvalues of $\widetilde F_p$ acting on $H_{\text{\'et}}^{2i}(\fX_{\overline{\Q}_p},\Q_{\ell})$ are $p^i$, for $s/2 \leq i \leq d$. The rest of the theorem follows by applying the $\ell$-adic \'etale to singular comparison theorem and the universal coefficients theorem. 
\end{proof}

Applying Proposition~\ref{prop:approximate} with $\fX = \ocM_{4,n}$ and $s = 9 + 2n$, we have the following.

\begin{corollary} \label{cor:approximate} 
Suppose there is a polynomial $P \in \Q[t]$ such that $\#\ocM_{4,n}(\F_q) = P(q) + o(q^{\frac{9}{2} + n})$, for $q = p^m$, as $m \to \infty$. Then $H^k(\ocM_{g,n}) = 0$ for all odd $k \leq 9$.
\end{corollary}

\noindent Similarly, to prove Theorem~\ref{thm:polynomial}, it suffices to show that the stated polynomial counts are correct up to $o(q^{\frac{9+n}{2}})$; this follows by applying Proposition~\ref{prop:approximate} with $\fX = \ocM_{4,n}$ and $s = 9 + n$.

\section{The polynomial point count property for the boundary divisor} \label{sec:boundary}

As explained in the introduction, our approach to showing that $\ocM_{4,n}$ has polynomial point count depends in an essential way on knowing that the boundary divisor $\partial \cM_{4,n}$ has polynomial point count. This is established inductively, as follows.

\begin{notation} \label{not:order}
Say $(g',n') \prec (g,n)$ if $(g',n') <_{\mathrm{lex}} (g,n)$ and $2g' + n' \leq 2g + n$.
\end{notation}

\noindent Each stratum of $\partial \cM_{g,n}$ is a finite quotient $\big( \prod_i \cM_{g_i,n_i} \big)  /G$, where each $(g_i, n_i) \prec (g,n)$.

\medskip

Say that a DM stack $\fX$ over $\Z$ has \emph{polynomial point count} if $\# \fX(\F_q)$ is a polynomial in $q$.

\begin{proposition} \label{prop:inductive-polynomial}
Suppose $\ocM_{g',n'}$ 
has polynomial point count for all $(g', n') \prec (g,n)$. Then
\begin{enumerate}
\item Every eigenvalue of $F_p$ on $H^\bullet_{c, \text{\'et}}((\cM_{g',n'})_{\overline \F_p}, \Q_\ell)$, for $(g',n') \prec (g,n)$, is a power of $p$;
\item The open moduli space $\cM_{g',n'}$ has polynomial point count for all $(g', n') \prec (g,n)$;
\item The boundary divisor $\partial \cM_{g,n}$ has polynomial point count.
\end{enumerate}
\end{proposition}

\begin{proof}
Assume $\ocM_{g',n'}$ has polynomial point count for all $(g', n') \prec (g,n)$. Since $\ocM_{g',n'}$ is smooth and proper over $\Z$, it follows that the odd cohomology of $\ocM_{g',n'}$ vanishes with rational coefficients, and every eigenvalue of $F_p$ acting on $H^{2k}_{\text{\'et}}((\ocM_{g',n'})_{\overline \F_p}, \Q_\ell)$ is $p^k$.

To prove (1), consider the weight spectral sequence associated to the normal crossing compactification $\cM_{g',n'} \subset \ocM_{g',n'}$, as in \cite[\S 2.3]{PayneWillwacher21}. The boundary is stratified according the topological types of the stable curves, which are encoded by the marked dual graphs, as in \cite{ArbarelloCornalba98}.  Let $\cM_G$ be the locus of curves with dual graph $G$, and let $\widetilde \cM_{G}$ be the normalization of its closure.  Then
\[
\widetilde \cM_{G} \cong \Big( \prod_{v \in V(G)} \ocM_{g_v, n_v} \Big) / \Aut(G),
\]
where $g_v$ and $n_v$ denote the genus and valence, respectively, of a vertex $v$, and hence
\[
H^\bullet(\widetilde \cM_G) = \Big(\bigotimes_{v \in V(G)} H^\bullet(\ocM_{g_v, n_v})\Big)^{\Aut(G)}.
\]
The weight spectral sequence abuts to the compactly supported cohomology of $\cM_{g,n}$, collapses at $E_2$, and has $E_1$-page given by
\begin{equation} \label{eq:weightjk}
E_1^{j,k} = \bigoplus_{|E(G)| = j} H^k(\widetilde \cM_G). 
\end{equation}
If we take this spectral sequence in $\ell$-adic \'etale cohomology, after basechange to $\overline \F_p$, then it is a spectral sequence of Galois representations \cite[Example~3.5]{Petersen17}. The odd rows vanish, and every eigenvalue of $F_p$ acting on the row $E^{\bullet, 2k}_1$ is $p^{k}$.  This proves (1).

Moreover, the dimension of $E_1^{j,k}$ is independent of the prime $p$ and  every eigenvalue of Frobenius acting on $H^\bullet_{c, \text{\'et}}((\cM_{g',n'})_{\overline \F_p}, \Q_\ell)$ is an integer power of $p$. The Behrend--Grothendieck--Lefschetz trace formula \eqref{eq:trace} tells us that $\cM_{g',n'}$ has polynomial point count. This proves (2).

Finally, note that $\partial \cM_{g,n}$ is the disjoint union of the spaces $\cM_G = \prod_{v \in V(G)}  \cM_{g_v, n_v} / \Aut(G)$, where $(g_v , n_v) \prec (g , n)$ for all $v$. The arguments above show that $\cM_G$ has polynomial point count for each graph $G$. This proves (3).
\end{proof}

To determine the polynomials $\#\ocM_{4,n}(\F_q)$ and $\#\cM_{4,n}(\F_q)$ from the approximate polynomial points that we compute in \S\S\ref{sec:cone}--\ref{sec:split}, we must know the precise point count for the boundary $\partial \cM_{4,n}$.  The computation of the latter involves the symmetric group action on the cohomology of $\cM_{g',n'}$ for $(g',n') \prec (4,n)$, due to the $\Aut(G)$-invariants in \eqref{eq:weightjk}.  See \S \ref{sec:equiv-boundary}.

\section{Counting hyperelliptic curves with up to three marked points} \label{sec:hypers}
Let $\cH_{g,n} \subset \cM_{g,n}$ be the closed substack of $n$-marked hyperelliptic curves, for $g \geq 2$.  Here, we recall the computation of  $\# \mathcal H_{g,n}(\F_q)$ for $n \leq 3$. We follow the arguments from \cite{Bergstrom09}, which explains such computations more generally for $n \leq 7$.  For simplicity, we present proofs only in the case where $q$ is odd, and give a reference for the general case.  

Let $q$ denote an odd prime power.  The stack $\cH_g$ of hyperelliptic curves of genus $g$ over $\Spec \Z\big[\frac{1}{2}\big]$ is a global quotient of the space of homogeneous polynomials of degree $2g+2$ with nonvanishing discriminant by a connected linear algebraic group $G$ with $\#G(\F_q) = (q^2-1)(q^2-q)$ \cite[Corollary~4.7]{ArsieVistoli04}. The discriminant vanishes if and only if the polynomial is divisible by the square of a nonconstant polynomial.  We therefore dehomogenize and consider the set $P_g$ of squarefree polynomials of degree $2g + 2$ or $2g + 1$ in a single variable over $\F_q$.  Since any polynomial can be written uniquely as the product of a square-free polynomial and the square of a monic polynomial, we get
\[ 
(q-1)(q^{2g+2}+q^{2g+1})=\# P_g+\# P_{g-1} \cdot q + \ldots +  \# P_1 \cdot q^{g-1}+ \#P_0 \cdot q^{g}+ (q-1)q^{g+1}. 
\]
Note that $\#P_0 = q^2(q-1)$. Then a simple induction shows $\# P_g= (q-1)(q^{2g+2}-q^{2g})$ for $g \geq 1$, and
\begin{equation}\label{eq:Hg}
\# \mathcal H_{g}(\F_q)= \frac{\# P_g}{\#G(\F_q)}=q^{2g-1},
\end{equation}
for $g \geq 2$; cf. \cite[Proposition 7.1]{BrockGranville01}. 

\begin{proposition} \label{prop:Hgn}
The point counts $\#\cH_{g,n}(\F_q)$ for $g \geq 2$ and $n \leq 3$ are
\[
\# \cH_{g,1}(\F_q) = (q+1) q^{2g-1}, \quad \quad \# \cH_{g,2}(\F_q) = (q+2)q^{2g}-1, \quad \quad \# \cH_{g,3}(\F_q) = (q^2+3q-1)q^{2g}-3q.
\]
\end{proposition}

We give a short proof in the case where $q$ is odd. The same polynomial formulas hold in general; see \cite[Theorem 10.3]{Bergstrom09}. Knowing these point counts for odd $q$ is sufficient for the purpose of proving our main theorems, since \cite[Theorem~2.1]{vandenBogaartEdixhoven05} only requires point counting information over a set of primes of Dirichlet density 1.

\begin{lemma} \label{lem:moments}
Assume $q$ is odd.  Then
\[
\sum_{C \in [\cH_g(\F_q)]} \frac{(q+1-\#C(\F_q))^k}{\# \Aut(C)} = \left\{ \begin{array}{ll} 0 &\mbox{ if $k$ is odd;} \\
q^{2g}-1 & \mbox{ if $k = 2$.} \end{array} \right.
\]
\end{lemma}

\begin{proof}
Let $C_f$ be the closure in $\P(1,1,g+1)$ of the affine curve given by $y^2 = f(x)$, for $f \in P_g$.  Since $q$ is odd, every hyperelliptic curve is of this form.  Fix a nonsquare $t \in \F_q$.  Then the map $C_f \mapsto C_{tf}$ induces an involution on $[\cH_g(\F_q)]$, and $\#\Aut(C_f) = \# \Aut(C_{tf})$.  Moreover
\[
q+1-\#C_f(\F_q) = -(q+1-\#C_{tf}(\F_q)).
\]
This proves the lemma when $k$ is odd.

It remains to consider the case $k = 2$. Let $\chi$ be the quadratic character on $\F_q$, so $\chi(a)$ is $0$ if $a$ is zero, $1$ if $a$ is a nonzero square and $-1$ if $a$ is a nonsquare.  Then 
\begin{equation*}
\begin{aligned}
\sum_{C \in [\cH_g(\F_q)]} \frac{(q+1-\#C(\F_q))^2}{\# \Aut(C)} & =\frac{1}{\#G(\F_q)} \sum_{f \in P_g} \Big( \sum_{x \in \P^1(\F_q)} \chi\big(f(x)\big)\Big)^2 \\
& = \frac{1}{\#G(\F_q)}  \sum_{f \in P_g} \Big( \sum_{x \in \mathbb P^1(\F_q) } \chi\big(f(x)\big)^2 + \sum_{x\neq y \in \mathbb P^1(\F_q)} \chi\big(f(x)\big)\chi\big(f(y)\big)  \Big).
\end{aligned}
\end{equation*}
Note that 
$\sum_{f \in P_g} \chi(f(x))^2$ 
is the number of polynomials in $P_g$ that do not vanish at $x$.  This is independent of $x$, and evaluating at $x = \infty$ shows 
\begin{equation*}
\begin{aligned}
\sum_{f \in P_g} \sum_{x \in \mathbb P^1(\F_q) } \chi\big(f(x)\big)^2 & = (q+1)(q-1)(q^{2g+2} - q^{2g+1}) \\ & = q^{2g} \cdot \#G(\F_q).
\end{aligned}
\end{equation*}
Also, if we take the sum over \emph{all} polynomials of degree at most $d \geq 2$
\[
\sum_{\deg(f) \leq d} \sum_{x \neq y} \chi\big(f(x)\big) \chi(f(y))
\]
then the total is zero, by interpolation.  An inductive argument, using again the fact that every polynomial can be written uniquely as the product of a square free polynomial and the square of a monic polynomial, shows that
\[
\sum_{f \in P_g} \sum_{x \neq y} \chi\big(f(x)\big) \chi\big(f(y)\big) = - \# G(\F_q),
\]
and the lemma follows.
\end{proof}

\begin{proof}[Proof of Proposition~\ref{prop:Hgn}]
Assume $q$ is odd. Applying Lemma~\ref{lem:moments} for $k = 1$ and \eqref{eq:Hg}, we have
\begin{equation*} 
 \# \mathcal H_{g,1}(\F_q) = \sum_{C \in [\mathcal H_g(\F_q)]} \frac{\# C(\F_q) }{\#\Aut(C)} =(q+1)q^{2g-1}.
\end{equation*}

\noindent Applying Lemma~\ref{lem:moments} also for $k = 2, 3$, we find 
\begin{equation*}
 \begin{aligned}
 \# \mathcal H_{g,2}(\F_q) &=\sum_{C \in [\mathcal H_g(\F_q)]} \frac{\# C(\F_q) \, \bigl(\# C(\F_q)-1 \bigr) }{\#\Aut(C)}, \\
 & =(q^2+q) \, \# \mathcal H_{g}(\F_q)+ \sum_{C \in [\cH_g(\F_q)]} \frac{(q+1-\#C(\F_q))^2}{\# \Aut(C)}, \\
& =(q+2)q^{2g}-1,
\end{aligned}
\end{equation*}
and 
\begin{equation*} \label{eq:Hg3}
 \# \mathcal H_{g,3}(\F_q)=(q^3-q) \, \# \mathcal H_{g}(\F_q)+3q \cdot \sum_{C \in [\cH_g(\F_q)]} \frac{(q+1-\#C(\F_q))^2}{\# \Aut(C)} =(q^2+3q-1)q^{2g}-3q,
\end{equation*}
as required.
\end{proof}

\section{Counting non-hyperelliptic curves via the classification of quadrics} \label{sec:strata}

Recall that the image of the canonical embedding of a non-hyperelliptic curve $C$ of genus 4 is the complete intersection of a quadric and a cubic in $\P^3$. Conversely, every smooth $(2,3)$-complete intersection in $\P^3$ is a canonically embedded curve of genus $4$.  Thus the moduli stack of non-hyperelliptic curves $\cM_{4} \smallsetminus \mathcal H_{4}$ is naturally identified with the quotient of the space of smooth $(2,3)$-complete intersections in $\P^3$ by the action of $\Aut(\P^3) = \PGL_4$.  Since $\PGL_4$ is a connected linear algebraic group, we can therefore count points over $\F_q$ in this moduli stack using Proposition~\ref{prop:pointcount}(ii), as follows. Let 
\[
S_{23}(\F_q) := \{ \mbox{smooth $(2,3)$-complete intersections in $\P^3$ over $\F_q$} \}.
\]
Then
\begin{equation} \label{eq:M4-H4}
\# (\cM_{4} \smallsetminus \mathcal H_4)(\F_q) \ = \ \# S_{23}(\F_q) \, / \, \# \PGL_4(\F_q).
\end{equation}
We compute the right hand side of \eqref{eq:M4-H4} by summing over cases according to the $\PGL_4(\F_q)$-orbit of the unique quadric containing the canonically embedded curve and using the orbit stabilizer theorem.  
 
The quadric containing a canonically embedded genus $4$ curve is necessarily reduced and irreducible.  There are precisely three $\PGL_4(\F_q)$-orbits of reduced and irreducible quadric surfaces in $\P^3$ over $\F_q$, represented by:
\begin{itemize}
\item a quadric cone $\Qcone \cong \P(1,1,2)$,  
\item a non-split smooth quadric $\Qns$ of Picard rank 1 over $\F_q$,
\item a split smooth quadric $\Qsplit \cong \P^1 \times \P^1$.
\end{itemize}
Note that $\Qns$ becomes isomorphic to $\P^1 \times \P^1$ after base change to $\F_{q^2}$. The nontrivial element of $\Gal(\F_{q^2} | \F_q)$ interchanges the two $\P^1$ factors; the geometric Frobenius is given by $$([x_0 \! : \! x_1], [y_0 \! : \! y_1]) \mapsto ([y_0^q \! : \! y_1^q], [x_0^q \! : \! x_1^q]).$$

\noindent Let $\Ncone(q)$ be the number of smooth curves of degree 6 in $\P(1,1,2)$ over $\F_q$.  Similarly, let $\Nns(q)$ and $\Nsplit(q)$ be the number of smooth curves in $\Qns$ and $\Qsplit$, resp., defined over $\F_q$ that that have geometric bidegree $(3,3)$ in $\P^1 \times \P^1$.

\begin{proposition} \label{prop:sumoverquadrics}
The number of non-hyperelliptic curves of genus $4$ over $\F_q$ is
\[
\# (\cM_{4} \smallsetminus \mathcal H_4)(\F_q) \ = \ \frac{\Ncone(q)}{\#\Aut(\Qcone)} + \frac{\Nns(q)}{\#\Aut(\Qns)} + \frac{\Nsplit(q)}{\#\Aut(\Qsplit)}.
\]
\end{proposition}

\begin{proof}
First, we note that each automorphism of $\Qcone$ extends to a linear automorphism of $\P^3$, because the linear series of hyperplane sections is complete.  Thus, $\Aut(\Qcone)$ is naturally identified with the subgroup of $\PGL_4(\F_q)$ that stabilizes $\Qcone$.  Furthermore, since the linear series of cubic sections is complete, $\Ncone(q)$ is the number of smooth (2,3)-complete intersections that lie on a fixed $\Qcone$.  Thus, by the orbit-stabilizer theorem, we have
\[
\frac{\Ncone(q)}{\#\Aut(\Qcone)} \ = \ \frac{\#\{C \in S_{23}(\F_q) : C \mbox{ lies on a quadric cone} \}}{\#\PGL_4(\F_q)}.
\]
Similarly,
\[
\frac{\Nns(q)}{\#\Aut(\Qns)} \ = \ \frac{\#\{C \in S_{23}(\F_q) : C \mbox{ lies on a non-split quadric}\}}{\#\PGL_4(\F_q)},
\]
and
\[ 
\frac{\Nsplit(q)}{\#\Aut(\Qsplit)} \ = \ \frac{\#\{C \in S_{23}(\F_q) : C \mbox{ lies on a split quadric}\}}{\#\PGL_4(\F_q)}.
\]
Hence, the proposition follows from \eqref{eq:M4-H4}, since every $C \in S_{23}(\F_q)$ lies on a unique quadric.
\end{proof}

Proposition~\ref{prop:sumoverquadrics} extends naturally to canonically embedded genus $4$ curves with $n$ marked points, for $n \leq 3$. Let $\Ncone_n(q)$ be the number of tuples $(C; p_1, \ldots, p_n)$ where $C$ is a smooth $(2,3)$-complete intersection on $\Qcone(q)$ and $p_1, \ldots, p_n$ are distinct $\F_q$-rational points on $C$, and similarly for $\Nns_n(q)$ and $\Nsplit_n(q)$.

\begin{proposition} \label{prop:pointedsumoverquadrics}
The number of $n$-pointed non-hyperelliptic curves of genus $4$ over $\F_q$ is
\[
\# (\cM_{4,n} \smallsetminus \mathcal H_{4,n})(\F_q) \ = \ \frac{\Ncone_n(q)}{\#\Aut(\Qcone)} + \frac{\Nns_n(q)}{\#\Aut(\Qns)} + \frac{\Nsplit_n(q)}{\#\Aut(\Qsplit)}.
\]
\end{proposition}

The orders of these automorphism groups are readily computed: 
\begin{equation} \label{eq:auts}
\#\Aut(\Qcone) = q^7 -q^6-q^5 + q^4; \quad \#\Aut(\Qns)= 2(q^6-q^2); \quad \#\Aut(\Qsplit) = 2(q^3 - q)^2.
\end{equation}

\noindent The following sections are devoted to estimating $\Ncone_n(q)$, $\Nns_n(q)$, and $\Nsplit_n(q)$ for $n \leq 3$, using a sieve for counting smooth curves in a family that allows for efficient computation of approximate point counts via properties of Hasse--Weil zeta functions.

\section{A Hasse--Weil sieve for counting smooth curves in families} \label{sec:HWSieve}

We begin by recalling the sieve method introduced for counting smooth curves in linear series over finite fields in  \cite[\S 6]{Bergstrom08}.  The sieve works equally well to count smooth fibers in arbitrary families of curves, not just linear series, and can be adapted in the evident way for families of higher dimensional varieties. Here, we focus on families of curves.

\smallskip

Let $\cC \to V$ be a surjective family of curves in a variety $X$ over $\F_q$.  We wish to count the smooth curves in this family, i.e., we want to determine $$N_\cC := \#\{v \in V(\F_q) : \cC_v \mbox{ is smooth}\}.$$ We do so by starting with $\#V(\F_q)$ and then adding and subtracting terms that account for singularities.  For a nonempty partition $\lambda$, let $X(\lambda)$ be the set of $0$-dimensional subsets  $Z \subset X$ such that the geometric Frobenius $F_q$ acts on $Z(\overline \F_q)$ with orbit type $\lambda$.  These are exactly the $\lambda$-tuples from \cite[Definition~4.7]{Bergstrom08}

\begin{proposition}[{\cite[\S 6]{Bergstrom08}}] \label{prop:naivesieve}
Suppose all fibers of $\cC \to V$ are reduced.  Then
\begin{equation} \label{eq:naivesieve}
N_{\cC} \, =  \, \#V(\F_q)  + \sum_\lambda \sum_{Z \in X(\lambda)}  (-1)^{\ell(\lambda)} \cdot \# \{v \in V(\F_q) : \cC_v \mbox{ is singular along $Z$}\}.
\end{equation}
\end{proposition}

\noindent The proof is elementary; if $C$ is a singular fiber then $\sum_{Z \subset C^\sing} (-1)^{\ell(\lambda_Z)} = -1$, where $Z$ ranges over nonempty subsets of $C^\sing$ defined over $\F_q$ and $\lambda_Z$ is the orbit type of $F_q$ acting on $Z(\overline \F_q)$.  

\begin{remark}
For applications to families that include non-reduced fibers, one approach is to modify the sieve and account for the non-reduced fibers separately, as in \cite[Definition~6.3]{Bergstrom08}. Under favorable circumstances, when the nonreduced locus of each fiber has only even cohomology, one may also apply the sieve directly and the contributions from nonreduced fibers eventually cancel.  See Proposition~\ref{prop:vanishsd}.
\end{remark}

Here, we reinterpret Proposition~\ref{prop:naivesieve} in terms of Hasse--Weil zeta functions to see that it can be extended to many families with non-reduced fibers, including the complete linear series on quadrics spanned by canonically embedded genus $4$ curves discussed in \S\ref{sec:strata}. We also discuss truncations of the sieve, obtained by summing over partitions of bounded length, and systematically controlling the error terms in the resulting approximations of $N_\cC$.  The approximations obtained in this way are typically better than one might naively expect. Indeed, as is well-known to practitioners of point-counting, there are often remarkable cancellations among the terms in sums such as \eqref{eq:naivesieve}. For instance, if we fix $d \geq 2$ and restrict the sum to partitions $\lambda \vdash d$, then the combined contribution of all partitions of $d$ is typically orders of magnitude smaller than the contribution of any given partition. These cancellations can be explained and systematically quantified by reinterpreting \eqref{eq:naivesieve} in terms of coefficients of inverse Hasse--Weil zeta functions and using fundamental facts about eigenvalues of Frobenius acting on compactly supported $\ell$-adic \'etale cohomology, as we now discuss.

\bigskip

Recall that the Hasse--Weil zeta function of a variety $Y$ over $\F_q$ is
\[
Z(Y;t) := \exp \Big( \sum_m \frac{\# Y(\F_{q^m})}{m}t^m \Big).
\]
Then an elementary argument, given by Vakil and Wood in the context of motivic discriminants, shows that the inverse of the Hasse--Weil zeta function satisfies
\begin{equation} \label{eq:HWZ2}
\frac{1}{Z(Y;t)} = \sum_d \Big( \sum_{\lambda \vdash d} (-1)^{\ell(\lambda)} \cdot \# Y(\lambda) \Big) \cdot t^d,
\end{equation}
where $\#Y(\lambda)$ denotes the number of Frobenius-stable subsets of $Y(\overline \F_q)$ with orbit type $\lambda$; see \cite[Proposition~3.7]{VakilWood15}.  The usefulness of this formula for point counting sieves and curve counting over finite fields was noted recently by Wennink \cite[p.~37]{Wennink20}.

\begin{definition}
Let $$s_d(Y) := \sum_{\lambda \vdash d} (-1)^{\ell(\lambda)} \cdot \# Y(\lambda)$$ be the coefficient of $t^d$ in the inverse Hasse--Weil zeta function of $Y$.
\end{definition}

Note that $s_0(Y) = 1$ for all $Y$.  Also, Hasse--Weil zeta functions and their inverses are multiplicative for disjoint unions, so $Z(\emptyset;t) = 1$, and
\begin{equation} \label{eq:HWmultiplicative}
\frac{1}{Z(Y;t)} = \frac{1}{Z(Y_1;t) \cdot Z(Y_2;t)}, \quad \mbox{ for } Y = Y_1 \sqcup Y_2.
\end{equation}
In particular, if $y \in Y$ is a point then 
\begin{equation} \label{eq:HWminuspoint}
\frac{1}{Z(Y \smallsetminus y; t)} = \frac{1}{1-t} \cdot \frac{1}{Z(Y; t)}.
\end{equation}
   
The Hasse--Weil zeta function is characterized in terms of the action of the geometric Frobenius $F_q$ on even and odd compactly supported \'etale cohomology, as follows
\begin{equation} \label{eq:HWFrob}
\frac{1}{Z(Y;t)} = \frac{\det\big( 1-t \, F_q \, | \, H^{\mathrm{even}}_{c,\text{\'et}}(Y_{\overline \F_q}, \Q_\ell) \big)}{\det\big( 1-t \, F_q \, | \, H^{\mathrm{odd}}_{c,\text{\'et}}(Y_{\overline \F_q}, \Q_\ell) \big)}.
\end{equation}
This is a consequence of the Grothendieck--Lefschetz trace formula \eqref{eq:trace}.

\begin{proposition} \label{prop:vanishsd}
Let $Y$ be a nonempty algebraic variety (possibly reducible or disconnected) that is proper over $\F_q$ and such that $H^{\mathrm{odd}}_{\text{\'et}}(Y_{\overline \F_q}, \Q_\ell) = 0$.  Then $$s_d(Y) = 0 \mbox{ \  for \ }  d > \dim H^{\mathrm{even}}_{\text{\'et}}(Y_{\overline \F_q}, \Q_\ell), \quad  \mbox{and} \quad  \sum_d s_d(Y) = 0.$$
\end{proposition}

\begin{proof}
By \eqref{eq:HWFrob}, the inverse Hasse--Weil zeta function $\frac{1}{Z(Y;t)}$ is a polynomial in $t$ of degree equal to $\dim H^{\mathrm{even}}_{\text{\'et}}(Y_{\overline \F_q}, \Q_\ell)$. Furthermore, it is divisible by $(1-t)$, since $1$ is an eigenvalue of $F_q$ acting on $H^0_{\text{\'et}}(Y_{\overline \F_q}, \Q_\ell)$.
\end{proof}

Recall that, when $X$ is an algebraic variety over $\F_q$, we write $X(\lambda)$ for the set of $F_q$-stable subsets of $X(\overline \F_q)$ with orbit type $\lambda$.  The following proposition, which we call the \emph{Hasse--Weil sieve} allows one to precisely count smooth fibers in a family of curves in $X$, even when some fibers are non-reduced, provided that the singular locus of each fiber has no odd cohomology.

\begin{proposition} \label{prop:HWsieve}
Let $\cC \to V$ be a surjective family of curves in a proper variety $X$ over $\F_q$.  Denote the number of smooth fibers by
\[
N_\cC = \# \{ v \in V(\F_q) : \cC_v \mbox{ is smooth}\}, 
\]
and let 
\begin{equation} \label{eq:bigsievecoeff}
S_d := \sum_{\lambda \vdash d} \sum_{Z \in X(\lambda)} (-1)^{\ell(\lambda)} \cdot \# \{ v \in V(\F_q) : \cC_v \mbox{ is singular along } Z\}.
\end{equation}  
Suppose $H^{\mathrm{odd}}_{\text{\'et}}((\cC^\sing_v)_{\overline \F_q}, \Q_\ell) = 0$, for all $v \in V(\F_q)$.  Then
\[
S_d = 0 \mbox{ for } d \gg 1  \quad \mbox{ and } \quad N_\cC = \#V(\F_q) + \sum_{d \geq 1} S_d.
\]
\end{proposition}

\noindent Note that the vanishing condition $H^{\mathrm{odd}}_{\text{\'et}}((\cC^\sing_v)_{\overline \F_q}, \Q_\ell) = 0$ is satisfied whenever $\cC_v$ is reduced, since then the singular locus is just a finite set of points.   

\begin{proof}
The right hand side of \eqref{eq:bigsievecoeff} may be rewritten as $\sum_{v \in V(\F_q)} s_d(\cC_v^\sing)$. By Proposition~\ref{prop:vanishsd}, we know that $s_d(\cC_v^\sing) = 0$ for $d$ sufficiently large. Moreover,
\[
\sum_{d \geq 0} s_d(\cC_v^\sing) = \left\{ \begin{array}{ll} 1 &\mbox{ if } \cC_v \mbox{ is smooth}; \\
0 & \mbox{ otherwise,} \end{array} \right.
\]
and the proposition follows.
\end{proof}

\begin{remark}
Proposition~\ref{prop:HWsieve} is useful not only for exactly counting the number of smooth fibers of $\cC \to V$, but also for efficiently approximating $N_\cC$.  Consider the truncation
\[
N_\cC[k] := \#V(\F_q) + \sum_{d \leq k} S_d.
\]
When the hypotheses of Proposition~\ref{prop:HWsieve} are satisfied, we have
\[
N_\cC - N_\cC[k] = \sum_{d > k} S_d.
\]
Then the terms $S_d$ for $d > k$ can be estimated by stratifying $V$ according to the topological type of the fibers.  The eigenvalues of $F_q$ acting on $H^{\mathrm{0}}_{\text{\'et}}((\cC^\sing_v)_{\overline \F_q}, \Q_\ell)$ are roots of unity, and the eigenvalues of $F_q$ acting on $H^{\mathrm{2}}_{\text{\'et}}((\cC^\sing_v)_{\overline \F_q}, \Q_\ell)$ are roots of unity times $q$.
\end{remark}

We conclude with examples showing that the hypothesis $H^{\mathrm{odd}}_{\text{\'et}}((\cC^\sing_v)_{\overline \F_q}, \Q_\ell) = 0$ is satisfied for the families of curves on quadrics that we will use for counting points in $(\cM_{4,n} \smallsetminus \mathcal H_{4,n})(\F_q)$, as discussed in \S\ref{sec:strata}.

\begin{example} \label{ex:nrcone}
Let $V \cong \A^{15}$ be the space of complete intersections of a quadric cone $\Qcone \subset \P^3$ with a cubic surface that does not contain the cone point, with $\cC \to V$ the restriction of the linear series to this open subspace.

The singular locus of any non-reduced fiber is isomorphic to $\P^1$, and the \'etale cohomology of a non-reduced scheme agrees with that of its reduced induced subscheme.  We compute that 
\[
s_0(\P^1) = 1; \quad s_1(\P^1) = -1 -q; \quad  s_2(\P^1) = q; \quad \mbox{ and } \quad s_d(\P^1) = 0 \mbox{ for } d > 2.
\]
In particular, Proposition~\ref{prop:HWsieve} applies, so we can count smooth fibers by applying the Hasse--Weil sieve to the family over $V$. Moreover, the contribution of each non-reduced curve to the truncated count $N_\cC[k]$ vanishes for $k \geq 2$. 

Let $\cC^\red \to V^\red$ be the restriction of $\cC \to V$ to the open locus $V^\red$ of reduced curves.  Thus, 
\[
N_\cC[k] = N_{\cC^{\red}}[k] \mbox{ for } k \geq 2.
\]
In particular, there is no need to modify the naive sieve given by Proposition~\ref{prop:naivesieve} to count smooth fibers in $\cC \to V$, as long as one groups terms appropriately, as in the Hasse--Weil sieve, and the estimates that we get from truncated counts are equally accurate, regardless of whether or not we include non-reduced fibers.
\end{example}

\begin{example} \label{ex:nrns}
Let $\Qns \subset \P^3$ be a non-split quadric over $\F_q$, and let $\cC \to V \cong \P^{15}$ be the linear series of complete intersections of $\Qns$ with cubics. The singular locus of any non-reduced fiber is either a curve of geometric bidegree $(1,1)$, which may be smooth or singular, or the disjoint union of such a $(1,1)$-curve with a point (the singular locus of the residual $(1,1)$-curve).  In all four cases, the odd cohomology of the singular locus vanishes, and hence we may apply Proposition~\ref{prop:HWsieve} to count the smooth fibers.  Moreover, since the cohomology of the singular locus of each non-reduced fiber has dimension at most 4, we have
$
N_\cC[k] = N^\red_\cC[k] 
$
for all $k \geq 4$. 
\end{example}

\begin{example} \label{ex:nrsplit}
Let $\Qsplit \subset \P^3$ be a split quadric over $\F_q$. We again consider the linear series of complete intersections with cubics $\cC \to V \cong \P^{15}$. In this case, the non-reduced locus is either a line in one of the rulings or a $(1,1)$-curve, which may be smooth or singular.  In both cases, the residual reduced curve (of type $(1,3)$ or $(1,1)$, respectively) may be smooth or singular. In all cases, the odd cohomology of the singular locus of each non-reduced fiber vanishes, and hence we may apply Proposition~\ref{prop:HWsieve} to count the smooth fibers.  Moreover, since the cohomology of the singular locus of each non-reduced fiber has dimension at most 5, we have
$
N_\cC[k] = N^\red_\cC[k] 
$
for all $k \geq 5$. 
\end{example}

\section{Applying the Hasse--Weil sieve to count genus 4 curves on quadrics}

We now count non-hyperelliptic curves of genus 4, by counting canonically embedded curves on each of the three isomorphism classes of irreducible quadric surfaces over $\F_q$, 
as outlined in \S\ref{sec:strata}.  In each case, we apply the Hasse--Weil sieve presented in \S\ref{sec:HWSieve} to count smooth fibers in an appropriate family of curves on the given surface.

\subsection{Canonically embedded genus 4 curves on a quadric cone} \label{sec:cone}

Let $\Qcone \subset \P^3$ be a quadric cone over $\F_q$, with $\Cone \to V^{\mathrm{con}} \cong \A^{15}$ the space of complete intersection curves of $\Qcone$ with a cubic that does not pass through the vertex of the cone.

Note that the smooth locus $X := (\Qcone)^{\mathrm{sm}}$ can be decomposed as $X = \A^2 \sqcup \A^1$, so
\begin{equation} \label{eq:ZQns}
\frac{1}{Z(X;t)} = (1-qt)(1-q^2t) = 1-(q^2 + q)t + q^3 t^2.
\end{equation}

As in Proposition~\ref{prop:HWsieve}, set 
\[
S_d(\Cone) := \sum_{\lambda \vdash d} \sum_{Z \in X(\lambda)} (-1)^{\ell(\lambda)} \cdot \#\big\{ v \in V^{\mathrm{con}}(\F_q) : Z \subset (\Cone_v)^\sing\big\}.
\]

\begin{proposition} \label{prop:Cone0}
The Hasse--Weil sieve terms $S_d(\Cone)$ for $d \leq 3$ are:
\[
S_0(\Cone) = q^{15}; \quad  S_1(\Cone) = -q^{14} -q^{13}; \quad S_2(\Cone) = q^{12}; \quad S_3(\Cone) = 0.
\]
\end{proposition}

\begin{proof}
By rearranging terms, we have
\[
S_d(\Cone) = \sum_{v \in V^{\mathrm{con}}(\F_q)} s_d((\Cone_v)^\sing).
\]
Thus $S_0(\Cone) = \# \{ v \in V^{\mathrm{con}}(\F_q) \} = q^{15}$.  Being singular at any point in $X(\F_q)$ imposes three linear conditions, so $$S_1 (\Cone)= q^{12} \cdot s_1(X) = -q^{14} - q^{13}.$$

Any fiber singular at two geometric points of a line in $\Qcone$ must contain the line and hence the vertex.  Being singular at two non-collinear geometric points imposes 6 linear conditions.  Thus $S_2(\Cone)$ is $q^9$ times the Hasse--Weil count of non-collinear pairs of points in $\Qcone$:
\[
S_2(\Cone) = q^9 \cdot (s_2(X) - (q+1) s_2(\A^1)) = q^{12}.
\]

Singularities at 3 non-collinear points impose $9$ conditions. It follows that $S_3(\Cone)$ is equal to $q^6 s_3(X) = 0$ plus terms divisible by $s_2(\A^1) = 0$ or $s_3(\A^1) = 0$ for triples containing 2 or 3 collinear points.  Thus $S_3(\Cone) = 0$.
\end{proof}

Let $V^\mathrm{con}_n \subset V^\mathrm{con} \times (\Qcone)^n$ be the incidence variety parametrizing tuples $(C; p_1, \ldots, p_n)$ where $p_1, \ldots, p_n$ are distinct points on $C$, and let $\Cone_n \to V^{\mathrm{con}}_n$ be the universal $n$-pointed curve.

\begin{proposition} \label{prop:ConeSieve}
The Hasse--Weil sieve terms $S_d(\Cone_n)$ for $d,n \leq 3$ are: 
\[
\resizebox{.95\hsize}{!}{
$\begin{array}{|l||l|}
\hline
S_0(\Cone_1) & q^{16} + q^{15}
\Tstrut \\ \hline
S_1(\Cone_1) &  -q^{15} -3q^{14} - q^{13} + q^{12}
\Tstrut \\ \hline
S_2(\Cone_1) &  2q^{13} +q^{12} - q^{11}
\Tstrut \\ \hline
S_3(\Cone_1) & 0
\Tstrut \\ \hline
\end{array}
\quad 
\begin{array}{|l||l|}
\hline
S_0(\Cone_2) & q^{17} + 2q^{16} -q^{14}
\Tstrut \\ \hline
S_1(\Cone_2) & -q^{16} -5q^{15} - 3q^{14} + 3q^{13} + 2q^{12}
\Tstrut \\ \hline
S_2(\Cone_2) & 3q^{14} +4q^{13} -3q^{12} - 3q^{11} + q^{10}
\Tstrut \\ \hline
S_3(\Cone_2) & -q^{12} +q^{11}+q^{10}-q^9
\Tstrut \\ \hline
\end{array}$
}
\]
\[
\begin{array}{|l||l|}
\hline
S_0(\Cone_3) & q^{18}+3q^{17}-5q^{15} -q^{14} + 2q^{13}
\Tstrut \\ \hline
S_1(\Cone_3) & -q^{17} -7q^{16} - 6q^{15} + 12q^{14} + 7q^{13} - 5q^{12}
\Tstrut \\ \hline
S_2(\Cone_3) & 4q^{15} +9q^{14} -10q^{13} - 9q^{12} + 6q^{11}
\Tstrut \\ \hline
S_3(\Cone_3) & -3q^{13} +5q^{12}-5q^{11}+7q^{10} + 2q^9 -12q^8+6q^7
\Tstrut \\ \hline
\end{array}
\]
\end{proposition}

\begin{proof}
Start with $n = 1$.  A point in $X$ imposes one linear condition on curves in $\Cone$, so $S_0(\Cone_1) = q^{14} (q^2+q)$.  The computations of $S_d(\Cone_1)$ for $1 \leq d \leq 3$ are similar to those of $S_d(\Cone)$, taking into account whether or not the singularity is at the marked point $x$, and using \eqref{eq:HWminuspoint} to compute $\frac{1}{Z(X \smallsetminus x;t)}$. Being singular at a point imposes $3$ conditions except when the singularity coincides with the marked point $x$, where a singularity is only 2 conditions, so
$$S_1(\Cone_1) =  q^{11}(q^2+q) s_1(X \smallsetminus x) + q^{12}(q^2+q) s_1(x).$$
The formula for $S_1(\Cone_1)$ follows, since $s_1(X \smallsetminus x) = -q^2-q+1$ and $s_1(x) = -1$. Similarly, letting $L$ be the line through the marked point, and discarding terms divisible by $s_2(\A^1) = 0$, we have
$$S_2(\Cone_1) = q^8 (q^2+q) \cdot (s_2(X \smallsetminus x) - s_2(L \smallsetminus x)) - q^9 \cdot (q^2+q) \cdot s_1(X \smallsetminus L).$$
The formula for $S_2(\Cone_1)$ follows, since $s_2(X \smallsetminus x) = q^3-q^2-q+1$ and $s_2(L\smallsetminus x) = 1-q$. 

For $n = 2$, note that passing through two points imposes two linear conditions, so $$S_0(\Cone_2) = (q^2+q)(q^2+q-1)q^{13}.$$ The remaining computations are again similar to those for $n = 1$, taking into account how many of the singularities coincide with marked points. For instance, $S_2(\Cone_{2})$ is equal to $$(q^2+q)(q^2)\cdot \big(q^7(s_2(X \smallsetminus 2 \mathrm{pt}) - 2s_2(\G_m)) + 2q^8(q^2 -1) + q^9\big),$$ coming from configurations where the marked points are not collinear, plus
$$(q+1)(q)(q-1)(2q^8)(q^2)$$ from configurations where the marked points are collinear.
Here, $2s_2(\G_m) = 2-2q$ is the Hasse--Weil count of configurations of two points collinear with, but distinct from, a marked point, and $2(q^2-1)$ is the count of non-collinear pairs including exactly one marked point. 

To compute $S_3(\Cone_2)$, note that the three singularities lie on a unique smooth hyperplane section. For a fixed hyperplane that contains at most 1 of the marked points, one finds that the contribution is a sum of terms each divisible by $s_3(\P^1) = 0$ or $s_2(\A^1) = 0$.  It remains to consider the cases where the hyperplane section contains both marked points. In such cases, the marked points are necessarily not collinear.  There are $(q^2+q)q^2$ choices for 2 non-collinear marked points, and $q$ smooth hyperplane sections through each such pair.  When the hyperplane is fixed, the Hasse--Weil count is 
\begin{equation}\label{eq:cone2}
s_1(\G_m) \cdot q^6 - 2 s_2(\G_m) q^5 + s_3(\G_m) q^5,
\end{equation}
where the term involving $s_k(\G_m)$ is the Hasse--Weil count of configurations of 3 singularities in which $k$ singularities are not at the marked points. Multiplying \eqref{eq:cone2} by $q^5+q^4$ gives the formula for $S_3(\Cone_2)$, using the fact that $s_d(\G_m) = 1-q$, for $d \geq 1$.

The computations for $\Cone_3$ are again similar, only more bookkeeping is required to keep track of how many singularities coincide with marked points, how many marked points are coplanar with three singularities, and so on.  
\end{proof}

\begin{remark} \label{rem:approximate}
These computations are more precise than what is needed for our main results. To prove Theorems~\ref{thm:polynomial} and \ref{thm:openpolynomial} via the approximate point counting, as in  
Proposition~\ref{prop:approximate},
we only need to know $\#S_d(\Cone_3)$ up to $bq^{13} + o(q^{13})$ for some (undetermined) integer $b$. For Theorem~\ref{thm:oddvanishing}, an approximate count up to $O(q^{14})$ would suffice. Such estimates are easier to obtain than the computations given here; many cases can be readily discarded.  Similarly, our point counting computations in \S\S\ref{sec:non-split}--\ref{sec:split} for curves whose canonical embeddings lie on smooth quadrics are far more precise than required for our main results.  
\end{remark}

\subsection{Counting canonical genus 4 curves on a non-split quadric} \label{sec:non-split}

Let $\Qns$ be a non-split quadric over $\F_q$, with $\Cns \to V^{\mathrm{nsp}} \cong \P^{15}$ the linear series of complete intersections with cubics in $\P^3$.

Applying formula \eqref{eq:HWFrob}, we see that the inverse Hasse--Weil zeta function of $\Qns$ is
\[
\frac{1}{Z(\Qns;t)} = (1-t)(1-qt)(1+qt)(1-q^2t) \]
and hence
\[
s_1(\Qns) = -q^2 - 1; \quad s_2(\Qns) = 0; \quad s_3(\Qns) = q^4 + q^2; \quad s_4(\Qns) = -q^4.
\]
The Hasse--Weil zeta function of $\Qns$ minus finitely many points is computed similarly, via \eqref{eq:HWminuspoint}.

We use the notation $\p_n := q^n + q^{n-1} + \cdots + 1$. In other words, $\p_n = \# \P^n(\F_q)$.

\begin{proposition} \label{prop:ns0}
The Hasse--Weil sieve terms $S_d(\Cns)$ for $d \leq 3$ are: 
\[
 S_0(\Cns) = \p_{15}; \quad S_1(\Cns) = (-q^2-1)\p_{12}; \quad S_2(\Cns) = 0; \quad S_3(\Cns) = (q^4 + q^2) \p_{6}.
  \]
\end{proposition}

\begin{proof}
For $n \leq 3$, singularities at any $n$ points impose $3n$ linear conditions, and hence $$S_n(\Cns) = s_n(\Qns) \p_{15-3n}.$$  
\vskip -18 pt \end{proof}

\smallskip

Let $V^\mathrm{nsp}_n \subset V^\mathrm{nsp} \times (\Qns)^n$ be the incidence variety parametrizing tuples $(C; p_1, \ldots, p_n)$ where $p_1, \ldots, p_n$ are distinct points on $C$, and let $\Cns_n \to V^{\mathrm{nsp}}_n$ be the universal $n$-pointed curve. 

\begin{proposition} \label{prop:NonsplitSieve}
The Hasse--Weil sieve terms $S_d(\Cns_n)$ for $d,n \leq 3$ are: 
\[
\resizebox{.95\hsize}{!}{
$\begin{array}{|l||l|}
\hline
S_0(\Cns_1) & (q^2+1)\p_{14}
\Tstrut \\ \hline
S_1(\Cns_1) & -(q^4+q^2)\p_{11} -(q^2+1)\p_{12}
\Tstrut \\ \hline
S_2(\Cns_1) & q^{13}+q^{11}
\Tstrut \\ \hline
S_3(\Cns_1) & (q^6+q^4)\p_5 + (q^4+q^2)\p_6
\Tstrut \\ \hline
\end{array}
\quad 
\begin{array}{|l||l|}
\hline
S_0(\Cns_2) & (q^4 + q^2) \p_{13}
\Tstrut \\ \hline
S_1(\Cns_2) & (-q^6+q^2)\p_{10} -2(q^4+q^2)\p_{11}
\Tstrut \\ \hline
S_2(\Cns_2) & 2 q^{14}+q^{13}+q^{12}+q^{11}-q^{10}
\Tstrut \\ \hline
S_3(\Cns_2) & 3 q^{11}+3 q^{10}+5 q^9+5 q^8+3 q^7+3 q^6+q^5+q^4
\Tstrut \\ \hline
\end{array}$
}
\]
\[
\begin{array}{|l||l|}
\hline
S_0(\Cns_3) & (q^6 - q^2) \p_{12}
\Tstrut \\ \hline
S_1(\Cns_3) & (-q^8 + 2 q^6 + q^4 - 2 q^2) \p_9 -3(q^6-q^2)\p_{10}
\Tstrut \\ \hline
S_2(\Cns_3) & 3(q^{15} + q^{14} -q^{13} -q^{11} -q^{10} + q^9)
\Tstrut \\ \hline
S_3(\Cns_3) & -2q^{13} +6q^{12} +2q^{11}+2q^{10} + 3q^9 -5q^8-2q^7-2q^6-q^5-q^4
\Tstrut \\ \hline
\end{array}
\]
\end{proposition}

\begin{proof}
The calculations for $d \leq 2$ are similar to those for $S_d(\Cns)$, with additional cases according to how many of the singularities are at the marked points.  For $d = 3$, the 3 singularities lie in a unique hyperplane section. When the hyperplane section is smooth, we consider how many of the marked points lie on the hyperplane, since a second or third marked point in the hyperplane imposes no new conditions, and how many of the singularities coincide with marked points. These computations are similar to the computation of $S_3(\Cone_2)$, above. 

If the hyperplane section is singular, then it is a conjugate pair of lines meeting at a point $x \in \Cns(\F_q)$.  The three singularities must be $x$ plus a conjugate pair of points, one on each line. Being singular at two points on a line forces containment of the line, so the curves under consideration consist of the hyperplane section plus a residual quadric section passing through the 2 conjugate points.  Any set of up to 3 marked points outside the hyperplane section imposes independent conditions on the residual quadric section, and the only possible marked point in the hyperplane is $x$.  There are $q^2 + 1$ possibilities for the singular hyperplane section, $q^2$ possibilities for the conjugate pair of singular points, and $q^2$ possibilities for marked points outside the hyperplane.  Thus, considering two cases according to whether or not $x$ is a marked point, we find that the contribution to $S_3(\Cns_1)$ from curves with 3 singularities spanning a singular hyperplane section is $q^2(q^2+1)(\p_6 + q^2 \p_5)$.  Similarly, the contributions to $S_3(\Cns_2)$ and $S_3(\Cns_3)$ are $q^2 (q^2 + 1)(q^2 \p_5 +  q^2 (q^2 - 1) \p_4)$ and $q^2 (q^2+1)(q^2(q^2-1) \p_4 + q^2(q^2-1)(q^2-2)\p_3)$, respectively.
\end{proof}

\subsection{Counting canonical genus 4 curves on a split quadric} \label{sec:split}

Let $\Qsplit$ be a split quadric over $\F_q$, with $\Csplit \to V^{\mathrm{spl}} \cong \P^{15}$ the linear series of complete intersections with cubics in $\P^3$.  Applying the formula \eqref{eq:HWFrob}, we see that the inverse Hasse--Weil zeta function of $\Qsplit$ is
\[
\frac{1}{Z(\Qsplit;t)} = (1-t)(1-qt)(1-qt)(1-q^2t) 
\]
and hence
\[
s_1(\Qsplit) = -q^2 - 2q - 1; \ \ \ s_2(\Qsplit) = 2(q^3+q^2+q); \ \ \ s_3(\Qsplit) = -q^4-2q^3-q^2.
\]
The Hasse--Weil zeta function of $\Qsplit$ minus finitely many points is computed similarly, via \eqref{eq:HWminuspoint}.

\begin{proposition} \label{prop:split0}
The Hasse--Weil sieve terms $S_d(\Csplit)$ for $d \leq 3$ are given by $S_0(\Csplit) = \p_{15}$,
\[
 S_1(\Csplit) = (-q^2-2q-1) \p_{12}; \ \  S_2(\Csplit) = 2(q^{3} + q^{2} + q) \p_9; \ \  S_3(\Csplit) = (-q^4-2q^3-q^2) \p_6.
  \]
\end{proposition}

\begin{proof}
Singularities at $d \leq 2$ points impose $3d$ independent conditions, so $S_d(\Csplit) = s_d(\Qsplit) \cdot \p_{15-3d}$.  Singularities at 3 points impose $9$ conditions unless the points lie on a line of one of the rulings, in which case they impose only 7 conditions.  When the 3 singularities are collinear in this way, the line is uniquely determined, and the contribution of such configurations is divisible by $s_3(\P^1) = 0$.  It follows that also $S_3(\Csplit) = s_3(\Qsplit) \cdot \p_{6}$.
\end{proof}

Let $V^\mathrm{spl}_n \subset V^\mathrm{spl} \times (\Qsplit)^n$ be the incidence variety parametrizing tuples $(C; p_1, \ldots, p_n)$ where $p_1, \ldots, p_n$ are distinct points on $C$, and let $\Csplit_n \to V^\mathrm{spl}_n$ be the universal $n$-pointed curve.

\begin{proposition} \label{prop:SplitSieve}
The Hasse--Weil sieve terms $S_d(\Csplit_n)$ for $d,n \leq 3$ are given up to $o(q^6)$ by: 
\[
\begin{array}{|l||l|}
\hline
S_0(\Csplit_1) & (q+1)^2\p_{14}
\Tstrut \\ \hline
S_1(\Csplit_1) &  - (q+1)^2((q^2+2q)\p_{11} + \p_{12})
\Tstrut \\ \hline
S_2(\Csplit_1) & 3q^{13} + 12q^{12} + 21q^{11} + 24q^{10} + 24q^9+ 24q^8 +24q^7 + 24q^6 + \cdots
\Tstrut \\ \hline
S_3(\Csplit_1) & -3q^{11} -10q^{10} -15q^9 -16q^8 - 16q^7-16q^6 + \cdots
\Tstrut \\ \hline
\end{array}
\]
\[
\begin{array}{|l||l|}
\hline
S_0(\Csplit_2) & (q+1)^2(q^2+2q) \p_{13}
\Tstrut \\ \hline
S_1(\Csplit_2) & -(q+1)^2(q^2+2q) ((q^2+2q-1)\p_{10} + 2\p_{11})
\Tstrut \\ \hline
S_2(\Csplit_2) & 4q^{14} + 27q^{13} + 61q^{12} + 75q^{11} + 73q^{10} + 72q^9 + 72q^8 + 72q^7 + 70q^6 + \cdots 
\Tstrut \\ \hline
S_3(\Csplit_2) & -10q^{12} -31q^{11} -45q^{10} -49q^9 - 49q^8 -47q^7 -41q^6 + \cdots 
\Tstrut \\ \hline
\end{array}
\]
\[
\begin{array}{|l||l|}
\hline
S_0(\Csplit_3) & (q+1)^2(q^2+2q)(q^2+2q-1) \p_{12}
\Tstrut \\ \hline
S_1(\Csplit_3) & -q^{17} - 12q^{16} - 54q^{15} - 106q^{14} - 117q^{13} - 102q^{12} - \\ & - 96q^{11} - 96q^{10} - 96q^{9} - 96q^8 - 95q^7 - 87q^6 + \cdots 
\Tstrut \\ \hline
S_2(\Csplit_3) & 5q^{15} + 53q^{14} + 151q^{13} + 190q^{12} + 153q^{11} + \\ & + 141q^{10} + 147q^9 + 142q^8 + 128q^7 + 90q^6 + \cdots
\Tstrut \\ \hline
S_3(\Csplit_3) & -30q^{13} -96q^{12} -116q^{11} - 94q^{10} - 97q^9 -93q^8 -62q^7 -24q^6 + \cdots
\Tstrut \\ \hline
\end{array}
\]
\end{proposition}

\begin{proof}
For the terms $S_1(\Csplit_n)$, note that a singularity and $n \leq 2$ marked points impose $3+n$ independent conditions, except when the singularity is equal to a marked point.  When $n = 3$, there is one additional case to consider: if the singularity and all 3 marked points are collinear, then they impose 5 conditions instead of 6.

For $S_2(\Csplit_1)$, the two singularities and a marked point impose independent conditions unless one singularity is equal to the marked point or all three points are collinear, i.e., they lie on a ruling. When all three are collinear and distinct, the contribution is divisible by $s_2(\A^1) = 0$. So, we need only consider two cases, according to whether or not one of the singularities is at the marked point.  

For $S_2(\Csplit_n)$ with $n = 2, 3$, there are more cases to consider, which require substantial bookkeeping but present no significant difficulties.  For instance, with $n = 3$ and when all 5 points are distinct, we distinguish cases where both singularities lie on a horizontal ruling and one singularity and all three marked points lie on a vertical ruling.

For $S_3(\Csplit_1)$, the computation is similar to that of $S_3(\Csplit)$.  For $S_3(\Csplit_2)$, we consider cases according to which collections of singularities and marked points are collinear or coplanar, and how many singularities coincide with marked points. Many such cases give contributions divisible by $s_3(\P^1)$ or $s_d(\A^1)$ for $d \geq 2$, and hence vanish.  One case with nonvanishing contribution is of particular note: there are two singular (1,1) curves through both marked points. Vanishing on the (1,1) curve imposes 7 conditions. Adding the three singularities gives 10 total conditions, instead of the expected 11.  The Hasse--Weil count of such configurations (for each such singular (1,1) curve) is $-s_1(\A^1)^2 = -q^2$, coming from configurations with one singularity at the singular point of the $(1,1)$-curve, and one additional singularity on each of its components. Another case with nonvanishing contribution is when all three singularities and both marked points are coplanar. The computation of $S_3(\Csplit_3)$ is similar, but with additional cases where a singularity is collinear with all three marked points.
\end{proof}

\subsection{Controlling the contribution from curves with 4 or more singularities}

In this section, we prove that the computations from the preceding sections suffice to determine $\#\ocM_{4,n}(\F_q)$ and hence $\#\cM_{4,n}(\F_q)$, for $n \leq 3$.  We do so, roughly speaking, by showing that the Hasse--Weil sieve count of curves with 4 or more singularities in each of the families that we considered is too small to be relevant when applying \cite[Theorem~2.1]{vandenBogaartEdixhoven05}.

\begin{proposition} \label{prop:error}
Suppose $\bullet \in \{ \mathrm{con}, \mathrm{nsp}, \mathrm{spl} \}$.  The remaining Hasse--Weil sieve terms are bounded by
\[
\sum_{d \geq 4} \frac{S_d(\cC^\bullet_n)}{\#\Aut(Q^\bullet)} = o(q^{\frac{9+n}{2}})
\]
for $n \leq 2$.  For $n = 3$, there is an integer $B$ such that 
\[
\sum_{d \geq 4} \frac{S_d(\cC^\bullet_n)}{\#\Aut(Q^\bullet)} = Bq^6 + o(q^6).
\]
\end{proposition}

\begin{remark}
For $n = 3$, we will not directly compute the coefficient $B$ of $q^6$.  Instead, we will use the fact that $\#\cM_{4,3}(\F_q)$ is a polynomial plus $o(q^6)$ to conclude that $\#\ocM_{4,3}(\F_q)$ is a polynomial and to determine all of the coefficients other than that of $q^6$.  We then use an Euler characteristic computation to determine the coefficient of $q^6$. We have also verified computationally that the resulting polynomial is correct at $q = 2, 3$; see Remark~\ref{rem:computations}. 
\end{remark}

\begin{proof}
First, we show that the contributions from nonreduced curves satisfy the specified bounds.  Example~\ref{ex:nrcone} shows that each nonreduced curve $C$ in $\Cone_n$ has $s_d(C) = 0$ for $d \geq 4$.

From Example~\ref{ex:nrns}, we see that each nonreduced curve $C$ in $\Cns_n$ has $s_d(C) = O(q^2)$ for $d \geq 4$.  The nonreduced curves in $\Cns_n$ form a family of dimension $6 + n$ and hence the total contribution is of size
\[
\frac{O(q^{8+n})}{\#\Aut(\Qns)} = O(q^{2+n}).
\]

The argument for nonreduced curves in $\Csplit_n$ is similar, using Example~\ref{ex:nrsplit}.  The contribution for curves nonreduced along a $(1,1)$-curve is bounded exactly as for $\Cns_n$.  The curves whose nonreduced locus is a line in a ruling form an $8$-dimensional family.  Each such curve $C$ has $s_d(C) = O(q)$ and hence the total contribution to $S_d(\Csplit_n)$ is
\[
\frac{b q^{9+n} + o(q^{9+n})}{\#\Aut(\Qsplit)} = bq^{3+n} + o(q^{3+n}),
\]
for some integer $b$, as required.

\smallskip

It remains to bound the contributions from reduced curves with 4 or more singularities.  We do so for curves on $\Qsplit$; the arguments for curves on $\Qcone$ and $\Qns$ are analogous and simpler. 

Start with the curves on $\Qsplit$ that have 5 or more singularities. Such curves are reducible and come in four families, with irreducible decompositions generically of the following types:
\begin{enumerate}[(A)]
\item $(2,1) + (1,2)$;
\item $(2,2) + (1,1)$, where the $(2,2)$ curve is singular;
\item $(3,2) + (0,1)$, where the $(3,2)$ curve has 2 singularities;
\item $(2,2) + (1,0) + (0,1)$. 
\end{enumerate}
Each of these families is of dimension $10$ and the general curve in each family has precisely $5$ singularities.  The contribution to $S_d(\Csplit_n)$ from any family of reduced curves of dimension $9$ trivially satisfies the stated bounds. Thus we can ignore reduced curves with 6 or more singularities as well as other degenerations, such as families of curves of type (A)-(D) with 4 or fewer singularities. Such curves necessarily have non-nodal singularities; they arise, e.g., when two components are tangent.  We now address the contributions from the families (A)-(D), and show that each contributes a polynomial of degree at most $9 + n$ to $S_d(\Csplit_n)$.

\smallskip

(A): Fix a $(2,1)$-curve $C$, which we may assume irreducible and hence smooth and isomorphic to $\P^1$.  There is a unique $(1,2)$-curve through $5$ general points on $C$. For $n = 0$, we see that the contribution of such curves is divisible by $s_5(\P^1) = 0$.  For $n = 1$, it is a sum of terms divisible by $s_5(\P^1)$ or $s_4(\A^1)$, which also vanishes.  For $n = 2$ and $3$, the contributions are still polynomial, and we consider cases depending on how many singularities coincide with marked points, taking a sieve count for the moving singularities in open subsets of $\P^1$. There are at least two moving singularities in each case, and the Hasse--Weil sieve counts of the moving singularities on $\P^1$ minus finitely many points is a polynomial of degree at most 1 in $q$, so we conclude that the contribution to $S_d(\Csplit_n)$ is a polynomial of degree at most $9 + n$, as required.

\smallskip

(B): The arguments are analogous to (A), fixing the singular point of the $(2,2)$-component, fixing the $(1,1)$-component, and varying the $4$ intersection points.

\smallskip

(C): This case is also analogous to (A), fixing the 2 singularities of the $(3,2)$-component, fixing the line, and varying the 3 intersection points on the line.

\smallskip

(D): Fix the two lines, and let $P$ be their point of intersection.  Choose
$2$ points away from $P$ on either ruling. We are interested
in the smooth $(2,2)$-curves intersecting the rulings
in the $4$ given points. These are curves of genus $1$
embedded with the two linear series of degree 2 and rank 1 given by the points on the
rulings. So the open set in the $\P^4$ of smooth $(2,2)$-curves through these 4 points 
is isomorphic to $\cM_{1,4}/(C_2\times C_2)$. It is known
that $\cM_{1,n}$ has polynomial point count for $n$ between
$4$ and $7$; this gives the required polynomiality.
The required cancellations arise since we vary $2$ points
on $\A^1$ and $s_2(\A^1) = 0$. In fact, we have \emph{two} copies of $\A^1$ on which
we vary $2$ points. For $n\leq2$, a single $\A^1$ suffices:
any choice of $2$ marked points outside the rulings
imposes $2$ conditions. For $n=3$, things are different:
when all $3$ marked points lie on a new ruling through
one of the $4$ given points, then (and only then) the third marked point
does not impose a new condition. The point on the original rulings
through which the new ruling passes is defined over $k$,
and the cancellation associated to the corresponding
original ruling no longer applies. However, the contribution of such curves is still divisible by $s_2(\A^1)$, by varying the two points on the other ruling.

\smallskip

Thus we have shown the contributions from nonreduced curves and from curves with 5 or more singularities (and their degenerations) are of the required form.  It remains to control the contributions from reduced curves with precisely 4 singularities, ignoring those that are degenerations of curves with 5 or more singularities.

\smallskip

Say that $4$ points on $\Qsplit \cong \P^1 \times \P^1$ are in \emph{general
position} if no $2$ of them lie on a ruling and they
do not all lie on an irreducible $(1,1)$-curve. If an
irreducible $(3,3)$-curve has $4$ singularities, the
singularities must be in general position. The reducible
$(3,3)$-curves with $4$ singularities that we must consider are of two
types: 
\begin{enumerate}[(A)]
\setcounter{enumi}{4}
\item $(2,2)+(1,1)$, where the $4$ singularities lie
on an irreducible $(1,1)$-curve;
\item $(3,2)+(0,1)$, where 
the $(3,2)$-curve is irreducible with $1$ singularity. The other $3$ singularities lie on a ruling. (Note that when $2$ singularities lie on a ruling, the
ruling necessarily contains a third singularity.)
\end{enumerate}
We first deal with the reducible curves. 

\smallskip

(E): Fix the $(1,1)$-curve. For $n = 0$, the signed count $S_d(\Csplit_n)$ is divisible by $s_4(\P^1) = 0$, by varying the singularities on the $(1,1)$-curve.  To see that this cancellation persists for $n \leq 3$, one only needs to understand in which situations an additional marked point fails to impose a condition on $(3,3)$-curves of this type.  There are two such cases.  First, if the marked point is on the $(1,1)$-curve, then the contribution is still divisibly by $s_4(\P^1)$.  The other case is when $n = 3$ and all $3$ marked points lie on a ruling through one of the 4 singularities.  The contribution in this case is divisible by $s_3(\A^1) = 0$, by varying the other $3$ singularities.

\smallskip

(F): We fix the horizontal ruling of type $(0,1)$ and the
singularity outside the ruling.
We may assume the vertical ruling through this
singularity does not contain another singularity (the $(3,3)$-curve
would in general have $5$ singularities).
Varying the 3 remaining singularities shows that the resulting contribution is divisible by $s_3(\A^1) = 0$. It remains to consider when $n \leq 3$ marked points fail to impose $n$ independent conditions.  This happens exactly
when $3$ marked points lie on the horizontal ruling through
the outside singularity, or when at least $2$ marked
points lie on the vertical ruling through the outside
singularity, or when $3$ marked points lie
on the vertical ruling through one of the singularities
on the fixed $(0,1)$-curve. Only in the latter case the choice
of the singularities is affected. Two of the singularities
can still vary freely on $\A^1$, and so the contribution is still divisible by $s_2(\A^1) = 0$, and hence the cancellation persists.

\smallskip

It remains to consider the irreducible curves with $4$ singularities in general position.  We claim that the Hasse--Weil count of such $4$-tuples is $(q+1)^2 q^2 (q-1)^2$, i.e., $\#(\PGL_2 \times \PGL_2)(\F_q)$.  Note that $\PGL_2 \times \PGL_2$ acts freely on such $4$-tuples, so the total must be divisible by the order of this group.  Also, the Hasse--Weil count of $4$-tuples on a $(1,1)$ curve is a sum of terms each divisible by $s_4(\P^1)$, or $s_d(\A^1)$ for $d\geq2$, and hence this count is 0.

The Hasse--Weil count of all $4$-tuples of points on $\Qsplit$ is $s_4(\Qsplit) = q^4$.  From this, we may sieve to remove configurations with at least two points on a line in the ruling.  We do so by yet another sieve, removing tuples with at least two points on each line in some configuration of lines over $\F_q$, with signs according to the number of Frobenius orbits on the set of lines, just as in Proposition~\ref{prop:naivesieve}. We consider cases, according to the number of lines in each ruling that contain at least 2 of the 4 points, the number of points on each line, and how many of the 4 points lie on multiple lines.  All cases contribute $O(q^5)$ except for the case of 1 line in each ruling, intersecting at $P$.  Here, all terms contribute 0 except where the $4$-tuple consists of $P$, one more point on each ruling, and one additional point. There are $(q+1)^2$ such pairs of lines, and each contributes $s_1(P) \cdot s_1(\A^1)^2 \cdot s_1(\A^2) = q^4$. Thus the full Hasse--Weil count of $4$-tuples in general position is $q^6+ O(q^5)$ and since it is divisible by $(q+1)^2 q^2 (q-1)^2$, it must be equal to $(q+1)^2 q^2 (q-1)^2$, as claimed.  There is a $\P^3$ of $(3,3)$-curves singular at any such configuration of $4$ points, and hence the total contribution of irreducible curves to $S_4(\Csplit)$ is $bq^{9} + o(q^{9})$.

Similarly, the curves of this type through configurations of $n$ marked points that impose independent conditions contribute $bq^{9+n} + o(q^{9+n})$.  It remains to account for the configurations of $n \leq 3$ points that impose fewer than $n$ independent conditions.  There are four possibilities to consider.  The first is when
$n = 3$, and all $3$ marked points lie on a ruling through
a singularity. A simple parameter count shows that these contribute $bq^{12} + o(q^{12})$.
The second possibility is when $n \geq 2$ and the $(1,1)$-curve
through $3$ of the singularities contains at least
$2$ marked points. Without cancellations, one finds
a polynomial of degree at most $13$, but the signed
count of the $3$ singularities on the $(1,1)$ gives more
cancellations than necessary.
The third possibility is when $n = 3$ and all $4$ singularities
and all $3$ marked points lie on a $(2,1)$ or $(1,2)$.
The argument is analogous to the one for the
second possibility.
Finally, the fourth possibility is when at least $1$ marked
point equals a singularity. When at least $2$ marked points
equal singularities, things are easy. When $1$ marked point
equals a singularity, the signed count of the remaining
$3$ singularities on $\Qsplit \smallsetminus \{ \mathrm{pt} \}$ 
gives the required cancellations.
\end{proof}

\begin{remark} \label{rmk:alreadypolynomial}
The counts carried out above are enough to conclude that $\#\ocM_{4,n}(\F_q)$ and $\#\cM_{4,n}(\F_q)$  are polynomials in $q$, for $n \leq 3$.  Determining these polynomials requires the precise computation of the contribution from the boundary, which is worked out in the next section, using equivariant point count data for $(g',n') \prec (g,n)$.
\end{remark}

\begin{remark}
There are many options for organizing point counts and carrying through the sieves. We have presented one approach, using the coefficients of inverse Hasse--Weil zeta functions systematically.  One can also carry out equivalent computations by more naive and elementary methods, e.g., in many cases one can simply apply Proposition~\ref{prop:naivesieve} directly, counting configurations in each Frobenius orbit-type.  Similarly, we used the Hasse--Weil zeta function to take a signed count of all $4$-tuples of points on $\Qsplit$ and then sieved to remove configurations with two or more points on a line, to obtain a signed count of $4$-tuples of points in general position.  But many other approaches to the same computation are possible. For instance, one could naively and directly count configurations in general position in each orbit type under Frobenius to find that the signed count is 
\begin{align*}
& 1/24 (q+1)^2 q^2 (q-1)^2 (q-2) (q-3)
-1/4 (q^2-q)^2 (q+1)^2 q (q-1) \\
& {}+1/3 (q^3-q)^2 (q+1) q
-1/4 (q^4-q^2) (q^3-q) (q-1) \\
& {}+1/8 (q^2-q)^2 (q^2-q-2) (q^2-2 q-3)
= (q+1)^2 q^2 (q-1)^2.
\end{align*}
\end{remark}

\section{Equivariant point counts and contributions from the boundary} \label{sec:equiv-boundary}

As explained in the introduction, our point counting strategy for determining $\# \ocM_{4,n}(\F_q)$ and $\# \cM_{4,n}(\F_q)$ for $n \leq 3$ depends in an essential way on knowing the corresponding boundary point counts $\# \partial \cM_{4,n}(\F_q)$.  The boundary point counts are calculated using a formula of Getzler and Kapranov for the characters of modular operads \cite[Theorem~8.3]{GetzlerKapranov98} that involves the $\SS_{n'}$-equivariant Euler characteristic of $\cM_{g',n'}$ for $(g',n') \prec (g,n)$.  Under favorable circumstances, this equivariant Euler characteristic can be computed via equivariant point counting \cite{KisinLehrer02, Bergstrom08}. We now recall the basics of this equivariant point counting in the form needed for the proof of our main results.

\subsection{Twisted forms and the trace formula}
Let $X$ be a variety over $\F_q$ with $\sigma \in \Aut(X)$.  Then there is a unique \emph{twisted form} of $X$, denoted $X^\sigma$ with an isomorphism $X^\sigma_{\overline \F_q} \xrightarrow{\sim} X_{\overline \F_q}$ that identifies the geometric Frobenius action on $X^\sigma_{\overline \F_q}$ with the action of $\sigma F_q$ on $X_{ \overline \F_q}$. 
We have already seen an example of such a twisted form in \S\ref{sec:strata}. The non-split quadric $\Qns$ is the twisted form of $\P^1 \times \P^1$ associated to the involution $\sigma$ that interchanges the two factors.  For a discussion of twisted forms of moduli spaces of curves obtained by permuting marked points, see \cite{KisinLehrer02, FlorenceHoffmannReichstein21}.

Because the endomorphism $\sigma F_q$ is identified with the geometric Frobenius of the twist $X^\sigma$, its action on $\ell$-adic \'etale cohomology satisfies the Grothendieck--Lefschetz trace formula. In particular, the graded trace of $\sigma F_q$ acting on $H^\bullet_{\text{\'et}}(X_{\overline \F_q},  \Q_\ell)$ is equal to $\#X^\sigma(\F_q)$.

\subsection{Equivariant point counts}

Let $X$ be a variety over $\F_q$ with the action of a finite group $G$. Note that conjugate elements of $G$ induce isomorphic twisted forms of $X$, so $\#X^\sigma(\F_q)$ is a class function on $G$.

\begin{definition}
The $G$-equivariant point count $\#^GX(\F_q)$ is the element of the representation ring of $G$ associated to the class function $\sigma \mapsto \# X^\sigma(\F_q)$.
\end{definition}

Note that the ordinary point count $\#X(\F_q)$ is the special case where $G$ is trivial. As with ordinary point counts, we will be particularly interested in how $\#^GX(\F_q)$ varies with $q$.

\begin{remark}
Although $X$ and its $G$-action may be defined over $\Z$, the construction of the twisted form $X^\sigma$ depends fundamentally on $q$.  Even in fixed characteristic, the construction does not commute with finite extension of the base field. For instance, the non-split quadric over $\F_{q^2}$ is not the base change of the non-split quadric over $\F_q$.  In this paper, we will always fix $q$ before constructing a twisted form, and we will count its points only over $\F_q$.  Thus, even though we omit $q$ from the notation for the twisted form, the meaning of $\#X^\sigma(\F_q)$ and $\#^GX(\F_q)$ is unambiguous.
\end{remark}

\begin{proposition} \label{prop:equiv-poly}
Suppose $X$ is smooth and proper over $\Z$ and has polynomial point count. Assume the action of $G$ is defined over $\Z$. Then $\#^GX(\F_q)$ is a polynomial in $q$.  More precisely, there are polynomials $P_V \in \Z[t]$, indexed by the irreducible representations $V$ of $G$, such that
\[
\#^GX(\F_q) = \sum_V P_V(q) \cdot [V].
\]
Moreover, if we write $P_V = \sum_i P_{V,i} t^i$ then $P_{V,i}$ is the multiplicity of $V$ in $H^{2i}(X_\C, \Q)$.
\end{proposition}

\noindent In particular, when $X$ is smooth and proper over $\Z$ and has polynomial point count, then we can determine $H^{\bullet}(X_\C, \Q)$ as a graded representation of $G$ by counting points over finite fields.

\begin{proof}
This follows from the universal coefficient theorem and the fact that the isomorphisms produced by the standard comparison theorems for singular and $\ell$-adic cohomology are functorial and hence $G$-equivariant \cite[Proposition~1.2]{KisinLehrer02}.
\end{proof}

\begin{remark} \label{rem:PD}
Let $V^\vee$ denote the irreducible representation dual to $V$.  When $X$ is irreducible of relative dimension $d$ over $\Z$, equivariant Poincar\'e duality tells us that $P_{V^\vee} (t) = t^d P_V(t^{-1})$.  In particular, when $V$ is self-dual, as is the case for all irreducible representations of $\SS_n$, the coefficients of $P_V$ have the symmetry $P_{V,i} = P_{V, d-i}$.
\end{remark}

The argument goes through essentially without change when $X$ is the coarse space of a smooth and proper DM stack over $\Z$. It also applies equally well to a smooth and proper DM stack  $\fX$, provided that one can define the twisted forms $\fX^\sigma$ for $\sigma \in G$.  For moduli spaces of stable curves with their symmetric group actions, this can be done as follows.

The moduli space $\ocM_{g,n}$ can be realized $\SS_n$-equivariantly as the global quotient of a smooth and proper variety $X_{g,n}$ by the action of a finite group $H$ \cite{BoggiPikaart00}. For fixed $q$, we can then define the twisted form $\ocM_{g,n}^\sigma := [X_{g,n}^\sigma / H]$, for $\sigma \in \SS_n$.  Note that $\ocM_{g,n}^\sigma$ is a smooth and proper DM stack over $\F_q$ and represents the functor taking a scheme $S$ over $\F_q$ to the set of isomorphism classes of families of stable curves over $S$ with $n$ disjoint labeled geometric sections such that Frobenius acts on this set of sections via the permutation $\sigma$.

\medskip

The irreducible representations of $\SS_n$ correspond to partitions $\lambda \vdash n$, and we identify these with symmetric functions in the usual way, writing $s_\lambda$ for the Schur function that corresponds to the character of the irreducible representation $V_\lambda$.  We write $\mathrm{gr}_\bullet^W$ for the associated graded of the weight filtration on the cohomology of a variety or DM stack.

\begin{corollary} \label{cor:equiv-polynomial}
Suppose $\ocM_{g',n'}$ has polynomial point count for $(g',n') \preceq (g,n)$.  Then there are polynomials $P_\lambda \in \Z[t]$, for $\lambda \vdash n$, such that
\[
\#^{\SS_n} \cM_{g,n}(\F_q) = \sum_{\lambda} P_\lambda(q) s_\lambda.
\]
Moreover, the coefficient of $t^k$ in $P_\lambda$ is the multiplicity of $V_\lambda$ in the virtual representation $\sum_i (-1)^i \mathrm{gr}_{2k}^W H^i_{c} (\cM_{g,n}, \Q)$.
\end{corollary}

\begin{proof}
By Proposition~\ref{prop:equiv-poly}, we know that $\#^{\SS_{n'}} \ocM_{g',n'}$ is a polynomial in $q$, for $(g',n') \preceq (g,n)$, with coefficients corresponding to the multiplicities of $V_\lambda$ in the cohomology groups of $\ocM_{g',n'}$. The rest of the argument is similar to the proof of Proposition~\ref{prop:inductive-polynomial}, by inspection of the weight spectral sequence for the normal crossings compactification $\cM_{g,n} \subset \ocM_{g,n}$s; see \cite[Example~3.5]{Petersen17} and \cite[\S 2.3]{PayneWillwacher21}.
\end{proof}

\noindent Thus, when $\ocM_{g',n'}$ has polynomial point count for $(g',n') \preceq (g,n)$, we can determine the $\SS_n$-equivariant weight-graded Euler characteristic of $\cM_{g,n}$ by counting points over finite fields.

\subsection{The contribution from the boundary} \label{sec:boundarycounts}
We compute the polynomials $\# \partial \cM_{4,n}(\F_q)$ for $n \leq 3$, using the $\SS_{n'}$-equivariant point counts on $\cM_{g',n'}$, for $(g',n') \prec (4,n)$ in the sense of Notation~\ref{not:order}, i.e., for $2g' + n' \leq 11$, and the Getzler-Kapranov formula for characters of modular operads \cite[Theorem~8.13]{GetzlerKapranov98}.

\begin{proposition} \label{prop:boundary-count}
The boundary point counts $\# \partial \cM_{4,n}(\F_q)$ are given by 
\begin{align*}
\#\partial \cM_4(\F_q) &= 3q^8 + 12q^7 + 33q^6 + 50q^5 + 50q^4 + 32q^3 + 13q^2 + 4q + 1; \\
\#\partial \cM_{4,1} (\F_q) & = 4q^9 + 28q^8 + 94q^7 + 192q^6 + 240q^5 +191q^4 + 93q^3 + 31q^2 + 6q + 1; \\
\#\partial \cM_{4,2} (\F_q) & =  8q^{10} + 72q^9 + 321 q^8 + 842 q^7 + 1362q^6 + 1362q^5 + 838q^4 + 321q^3 +   
\\  & \quad 
{}+ 78q^2 + 11q + 1; \\
\#\partial \cM_{4,3} (\F_q) & =17q^{11} + 200q^{10} + 1172 q^9 + 3990 q^8 + 8292 q^7 + 10606q^6 + 8296q^5 +  
\\ & \quad 
{}+ 3977 q^4 + 1179q^3 + 205 q^2 + 19 q + 2.
\end{align*} 
\end{proposition}
\begin{proof} Fix any $q$. The Getzler-Kapranov formula \cite[Theorem~8.13]{GetzlerKapranov98} gives a way to compute the ($\SS_n$-equivariant) weight-graded Euler characteristic of $\partial \cM_{g,n}$ using the $\SS_{n'}$-equivariant weight-graded Euler characteristics of $\cM_{g',n'}$, for $(g',n') \prec (g,n)$. We follow their notation.  Using the Lefschetz trace formula, it is straightforward to generalize the Getzler-Kapranov formula to point counts, noting that $p_n \circ q=q^n$, where $p_n$ is the degree $n$ power sum symmetric function and $\circ$ denotes plethysm. Namely, if we put 
 \[
\mathbb{C}\mathrm{h}(\mathcal V)=\sum_{2(g-1)+n>0} h^{g-1} \, \#^{\SS_n} \cM_{g,n}(\F_q),
\]
with $h$ a formal variable, then we can use this information to compute 
\[
\mathbb{C}\mathrm{h}(\mathbb{M} \mathcal V)-\mathbb{C}\mathrm{h}(\mathcal V)=\sum_{2(g-1)+n>0} h^{g-1} \,  \#^{\SS_n}\partial \cM_{g,n}(\F_q).
\]
More precisely, knowing $\#^{\SS_n'}\ \cM_{g',n'}(\F_q)$ for all $(g',n') \prec (g,n)$ we can determine 
$\#^{\SS_n}\ \partial \cM_{g,n}(\F_q)$. The required input for $g' \leq 3$ is readily found in the literature; see \cite{KisinLehrer02} for $g' = 0$ and \cite{Getzler98, Bergstrom08, Bergstrom09} for $1 \leq g' \leq 3$ and $(g',n') \prec (4,3)$.

Thus, we can use the Getzler-Kapranov formula plus previously known results in lower genus to compute $\# \partial \cM_{4}(\F_q)$. These lower genus results also suffice to compute $\# \partial \cM_{4,1}(\F_q)$.  

For $n = 2$ and $3$, some additional input is required; there is no circularity, one proceeds inductively, increasing the number of marked points one step at a time.  We use $\# \partial \cM_{4,1}(\F_q)$ plus an approximation to $\#\cM_{4,1}(\F_q)$ to determine $\#\ocM_{4,1}(\F_q)$ via Poincar\'e duality; for details see \S\ref{sec:mainproof}, below.  Subtracting $\# \partial \cM_{4,1}(\F_q)$ then gives a precise formula for $\#\cM_{4,1}(\F_q)$, which is used as input to compute $\# \partial \cM_{4,2}(\F_q)$, and so on.  Note that the equivariant point count $\#^{\SS_2} \cM_{4,2}(\F_q)$ is not required to determine $\# \partial \cM_{4,3}(\F_q)$. The ordinary point count $\# \cM_{4,2}(\F_q)$ suffices because there is only one stable graph of genus $4$ with $3$ marked points that has a vertex of genus 4 and valence 2, and it does not have any nontrivial automorphisms.

The actual computations were carried out using the computer software Maple and the symmetric polynomial package SF \cite{SF}.
\end{proof}

\section{Proof of Theorems~\ref{thm:polynomial} and \ref{thm:openpolynomial}} \label{sec:mainproof}

In this section, we discuss how to put together the approximate point counts on $\cM_{4,n}(\F_q)$ with the boundary point counts from the previous section to obtain precise point counts for $\ocM_{4,n}$ and $\cM_{4,n}$, for $n \leq 3$ and thus prove our main results stated in the introduction.

\subsection{Proof of Theorems~\ref{thm:polynomial} and \ref{thm:openpolynomial}} 
Consider first the case $n = 0$.  By Propositions~\ref{prop:Cone0} and \ref{prop:error}, we have
\[
\Ncone(\F_q) = q^{15} -q^{14} - q^{13} + q^{12} + o(q^{\frac{23}{2}}).
\]
Similarly, Propositions~\ref{prop:ns0}, \ref{prop:split0}, and \ref{prop:error} give
\[
\Nns(\F_q) = q^{15} -q^{12} -q^{11} + o(q^{\frac{21}{2}})  \quad \quad  \mbox{ and } \quad \quad \Nsplit = q^{15} - 2q^{13} -q^{12} + q^{11} + o(q^{\frac{21}{2}}).
\]
By Proposition~\ref{prop:pointedsumoverquadrics} we can get the count of all non-hyperelliptic curves by dividing these counts by the orders of the respective automorphism groups, as given in \eqref{eq:auts}, and summing to get
\[
\# (\cM_{4} \smallsetminus \cH_4)(\F_q) = q^9 + q^8 - q^6  +o(q^{\frac{9}{2}}).
\]
We then add $\# \cH_4(\F_q) = q^7$ and $\# \partial \cM_4(\F_q) = 3q^8 + 12q^7 + 33q^6 + 50q^5 + o(q^{\frac{9}{2}})$ to find 
\[
\# \ocM_4(\F_q) = q^9 + 4q^8  + 13q^7 + 32q^6 + 50q^5 + o(q^{\frac{9}{2}}).
\]
By \cite[Theorem~2.1]{vandenBogaartEdixhoven05}, it follows that $\#\ocM_{4}(\F_q)$ is a polynomial in $q$ and, by Poincar\'e duality, the polynomial must be
\[
\# \ocM_4(\F_q) = q^9 + 4q^8  + 13q^7 + 32q^6 + 50q^5 + 50q^4 + 32q^3 + 13q^2 + 4q + 1.
\]
Subtracting $\# \partial \cM_4(\F_q)$ then gives
\[
\# \cM_4(\F_q) = q^9 + q^8 + q^7 - q^6,
\]
as required.  The computations for $n = 1$ and $n = 2$ are similar.  Note that, for $n = 2$, one uses $\#\cM_{4,1}(\F_q)$ as input for calculating $\# \partial \cM_{4,2}(\F_q)$.

For $n = 3$, the analogous computations show that
\[
\# \ocM_{4,3}(\F_q) = q^{12}+21q^{11}+207q^{10}+1168q^9+3977q^8+8296q^7+Bq^6 + o(q^6),
\]
for some integer $B$.  Again,  \cite[Theorem~2.1]{vandenBogaartEdixhoven05} shows that $\#\ocM_{4,3}(\F_q)$ is a polynomial in $q$. Applying Poincar\'e duality for $\ocM_{4,3}(\F_q)$ and subtracting $\# \partial \cM_{4,3}(\F_q)$ (which is determined using $\#\cM_{4,2}(\F_q)$) then gives
\[
\# \cM_{4,3}(\F_q) = q^{12} + 4q^{11} + 7q^{10} -4q^9 -13q^8 + 4q^7 +(B-10606)q^6 - 11q^4 + 2q^2 + 2q - 1.
\]
Evaluating at $q = 1$ gives the Euler characteristic, and $\chi(\cM_{4,3}) = -10$, e.g., by \cite{Gorsky14}. We conclude that $B = 10605$, and the theorems follow.

\section{$\SS_n$-equivariant point counts on $\cM_{4,n}$ and $\ocM_{4,n}$ for $n = 2, 3$} \label{sec:S_n-equiv-M4n}

As we have seen, knowing the $\SS_{n'}$-equivariant point counts on $\cM_{g',n'}$ for $(g',n') \prec (g,n)$ is essential for determining the boundary point count $\# \partial \cM_{g,n}(\F_q)$. Moreover, we know that the equivariant point counts on $\ocM_{4,2}$, $\ocM_{4,3}$, $\cM_{4,2}$ and $\cM_{4,3}$ are polynomial, by Proposition~\ref{prop:equiv-poly} and Corollary~\ref{cor:equiv-polynomial}.  For future work, it will be essential to know these polynomials. 

\begin{theorem} \label{thm:equiv-point-counts}
The $\SS_n$-equivariant point counts on $\ocM_{4,n}(\F_q)$ and $\cM_{4,n}(\F_q)$, for $n = 2, 3$ are 
\[
\resizebox{.95\hsize}{!}{
$
\begin{array}{|l||l|}
\hline
\mbox{\raisebox{-1.5pt}{$\#^{\SS_2}(\ocM_{4,2})(\F_q)$}} & \mbox{\raisebox{-1.5pt}{$(q^{11}+9\,q^{10}+55\,q^9+220\,q^8+561\,q^7+901\,q^6+\cdots 
+9\,q+1) \, s_{2} \,+$}}  \\ & 
\mbox{\raisebox{-.5pt}{$+  (2\,q^{10}+21\,q^9+99\,q^8+277\,q^7+461\,q^6+\cdots 
+21\,q^2+2\,q) \, s_{1^2}$}}
\Tstrut \\ \hline
\mbox{\raisebox{-1.5pt}{$\#^{\SS_3}(\ocM_{4,3})(\F_q)$}} &  \mbox{\raisebox{-1.5pt}{$(q^{12}+11\,q^{11}+87\,q^{10}+424\,q^9+1347\,q^8+2694\,q^7+3414\,q^6+\cdots 
+1)\,s_{3} \, +$}}\\
& \mbox{\raisebox{-1pt}{$+ (5\,q^{11}+58\,q^{10}+349\,q^9+1220\,q^8+2578\,q^7+3304\,q^6+\cdots 
+5\,q)\,s_{21} \, +$}}\\
& \mbox{\raisebox{-.5pt}{$+ (4\,q^{10}+46\,q^9+190\,q^8+446\,q^7+583\,q^6+\cdots 
+46\,q^3+4\,q^2)\,s_{1^3}$}}
\Tstrut \\ \hline
\end{array}$}
\]

\[
\begin{array}{|l||l|}
\hline
\#^{\SS_2}(\cM_{4,2})(\F_q) &  (q^{11}+2q^{10}+3q^9-2q^8-2q^7-q^3-q^2)s_2+(q^{10}+q^9-2q^7-q^3-q^2)s_{1^2} 
\Tstrut \\ \hline
\mbox{\raisebox{-1pt}{$\#^{\SS_3}(\cM_{4,3})(\F_q)$}} &  \mbox{\raisebox{-1pt}{$(q^{12}+2q^{11}+3q^{10}-2q^9-4q^8+2q^7-q^6-q^3+2q^2)s_3+$}} \\ & 
\mbox{\raisebox{-.5pt}{$+(q^{11}+2q^{10}-q^9-4q^8+q^7-4q^3-1)s_{2,1}+(-q^8-2q^3+2q+1)s_{1^3}$}}
\Tstrut \\ \hline
\end{array}
\]
\end{theorem}

\noindent The elided terms above are determined by Poincar\'e duality, as in Remark~\ref{rem:PD}. In principle, Theorem~\ref{thm:equiv-point-counts} could be proved by direct calculation, evaluating at several values of $q$. However, doing so by brute force is beyond the limitations of current computers. Instead, we give an argument that is similar to the proofs of Theorems~\ref{thm:polynomial} and \ref{thm:openpolynomial}.

Let $\lambda \vdash n$ be a partition.  Say $\lambda$ has $n_1$ parts of size $1$, $n_2$ parts of size $2$, etc., and let $\sigma(\lambda) \in \SS_n$ be the permutation that fixes the first $n_1$ elements of $\{1, \ldots, n\}$, transposes the next $n_2$ pairs, and so on. Let
\[
\ocM_{g,\lambda} := \ocM_{g,n}^{\sigma(\lambda)},
\] 
be the corresponding twisted form. Note that a point in $\ocM_{g,\lambda}(\F_q)$ is represented by a tuple $(C; p_1, \ldots, p_n)$, where $C$ is a smooth curve of genus $g$ over $\F_q$ and $p_1, \ldots, p_n$ are points in $C(\overline \F_q)$ on which the geometric Frobenius $F_q$ acts by fixing the first $n_1$ points, transposing the next $n_2$ pairs, and so on.  We define $\cM_{g, \lambda}$, $\cH_{g,\lambda}$, and $\partial \cM_{g,\lambda}$  similarly.  In particular, $\ocM_{g, \lambda} = \cM_{g, \lambda} \sqcup \partial \cM_{g, \lambda}$. 

We compute $\# \ocM_{g,\lambda}(\F_q)$, for odd $q$, just as we computed $\# \ocM_{g,n}(\F_q)$:
\begin{itemize}
\item The hyperelliptic point count $\# \cH_{g, \lambda}(\F_q)$ is previously known \cite{Bergstrom09}.
\item The boundary point count $\# \partial \cM_{g,\lambda}(\F_q)$ is determined from the $\SS_{n'}$-equivariant point counts on $\cM_{g',n'}$ for $(g',n') \prec (g,n)$, via \cite[Theorem~8.13]{GetzlerKapranov98}.\;
\item The count of smooth non-hyperelliptic curves is approximated to the required accuracy by applying the Hasse--Weil sieve to $\sigma(\lambda)$-twisted forms of the universal families of $n$-pointed complete intersections on each of the three isomorphism classes of reduced and irreducible quadrics over $\F_q$, and then dividing by the orders of the respective automorphism groups of these surfaces.
\end{itemize}
Combining these data with $\SS_n$-equivariant Poincar\'e duality for $\ocM_{g,n}$ plus the $\SS_n$-equivariant Euler characteristic of $\cM_{g,n}$, known from \cite{Gorsky14}, is enough to determine the $\SS_n$-equivariant point count on $\ocM_{g,n}$.  Finally, the $\SS_n$-equivariant point count on $\cM_{g,n}$ is obtained by subtracting the contribution from the boundary.

\subsection{Equivariant count of hyperelliptic curves} \label{sec:equivarianthypers}

In \S\ref{sec:hypers}, we explained the count of hyperelliptic curves $\#\cH_{g,n}(\F_q)$ for $n \leq 3$, following \cite{Bergstrom09}, where such counts are carried out in much greater generality, $\SS_n$-equivariantly for $n \leq 7$.  For $n = 2, 3$, the $\SS_n$-equivariant counts are:
\[
\begin{array}{|l||l|}
\hline
\#^{\SS_2} \mathcal H_{4,2}(\F_q)&  (q^9+q^8-1)\,s_{2}+q^8\,s_{1^2}
\Tstrut \\ \hline
\#^{\SS_3} \mathcal H_{4,3}(\F_q) &  (q^{10}+q^9-q^8-q)\,s_3+ (q^9-q)\,s_{2,1}
\Tstrut \\ \hline
\end{array}
\]

\subsection{Equivariant count of canonically embedded genus $4$ curves}
 
We give the $\SS_n$-equivariant count of $n$-pointed canonically embedded genus $4$ curves over $\F_q$ by separately counting complete intersections with cubics on each of the three isomorphism classes of irreducible quadrics over $\F_q$, just as we have done for the ordinary point counts in \S\S\ref{sec:cone}--\ref{sec:split}.

Continuing the partition notation established above, we fix $\lambda \vdash n$ with $n_1$ parts of size $1$, $n_2$ parts of size $2$, etc., and let $\sigma(\lambda) \in \SS_n$ be the permutation that fixes the first $n_1$ elements of $\{1, \ldots, n\}$, transposes the next $n_2$ pairs, and so on. Let
\[
V_\lambda^\mathrm{con} := (V^{\mathrm{con}}_n)^{\sigma(\lambda)},
\] 
be the corresponding twisted form of the base $V^\mathrm{con}_n$ of the family of $n$-pointed complete intersections of a cubic with the quadric cone $\Qcone$.  Thus $V_\lambda^\mathrm{con}(\F_q)$ is the set of tuples $(C; p_1, \ldots, p_n)$ with $p_i \in C(\overline \F_q)$ such that the geometric Frobenius fixes the first $n_1$ points, transposes the next $n_2$ pairs, and so on. Let
\[
\Cone_\lambda \to V_\lambda^\mathrm{con} 
\]
be the universal $n$-marked complete intersection of $\Qcone$ with a cubic not passing through the vertex of the cone.  We define $\Cns_\lambda \to V_\lambda^\mathrm{ns}$ and $\Csplit_\lambda \to V_\lambda^\mathrm{spl}$ similarly.

Just as in the nonequivariant counts, we estimate the number of smooth fibers $N^\bullet_\lambda$  in the family $\cC^\bullet_\lambda \to V^\bullet_\lambda$, by applying the Hasse--Weil sieve. The sieve terms satisfy the conclusion of Proposition~\ref{prop:error}, with essentially the same proof, and hence, using $\SS_n$-equivariant Poincar\'e duality and, for $n = 3$, Gorsky's computation of the $\SS_n$-equivariant Euler characteristic of $\cM_{g,n}$ \cite{Gorsky14}, it suffices to compute the sieve terms $S_d(\cC^\bullet_\lambda)$  for $\lambda \in \{ [2], [2,1], [3] \}$ and $d \leq 3$.

\begin{proposition} \label{prop:ConeEquiv}
The sieve terms $S_d(\Cone_\lambda)$ for $\lambda \in \{ [2], [2,1], [3] \}$ and $d \leq 3$ are 
\[
\begin{array}{|l||l|}
\hline
S_0(\Cone_{[2]}) &  q^{17} - q^{14}
\Tstrut \\ [2.5 pt] \hline
S_1(\Cone_{[2]}) &  -q^{16} -q^{15} + q^{14} + q^{13}
\Tstrut \\ [2.5 pt] \hline
S_2(\Cone_{[2]}) &  q^{14}-2q^{13} + q^{12} + q^{11} - q^{10}
\Tstrut \\ [2.5 pt] \hline
S_3(\Cone_{[2]}) &  q^{12} - q^{11} -q^{10} + q^9
\Tstrut \\ [2.5 pt] \hline
\end{array}
\ \ 
\begin{array}{|l||l|}
\hline
S_0(\Cone_{[2,1]}) &  q^{18} +q^{17} - q^{15} - q^{14}
\Tstrut \\ [2.5 pt] \hline
S_1(\Cone_{[2,1]}) &   - q^{17}-3 q^{16}+4 q^{14}+q^{13}-q^{12} 
\Tstrut \\ [2.5 pt] \hline
S_2(\Cone_{[2,1]}) &  2 q^{15}-q^{14}-4 q^{13}+3 q^{12}+2 q^{11}-2 q^{10} 
\Tstrut \\ [2.5 pt] \hline
S_3(\Cone_{[2,1]}) &  q^{13} + q^{12} - 5q^{11} + q^{10} + 4q^9 -2q^8
\Tstrut \\ [2.5 pt] \hline
\end{array}
\]
\[
\begin{array}{|l||l|}
\hline
S_0(\Cone_{[3]}) &  q^{18}+q^{15}-q^{14}-q^{13}
\Tstrut \\ [2.5 pt] \hline
S_1(\Cone_{[3]}) &  -q^{17}-q^{16} + q^{13} + q^{12}
\Tstrut \\ [2.5 pt] \hline
S_2(\Cone_{[3]}) &  q^{15} - q^{13}
\Tstrut \\ [2.5 pt] \hline
S_3(\Cone_{[3]}) &  -q^{12} + q^{11} + q^{10} -q^9
\Tstrut \\ [2.5 pt] \hline
\end{array}
 \ .
\]
\end{proposition}

\begin{proof}
The proof is similar to that of Proposition~\ref{prop:ConeSieve} for $n = 2$ and $3$. The computations here are somewhat easier because there are fewer cases to consider.  For instance, when $\lambda = [3]$, the number of marked points that coincide with singularities must be either 0 or 3. 

Put $X := (\Qcone)^{\mathrm{sm}}$. To compute $S_0(\Cone_{[2]})$, note that there are $q^4 - q$ ordered conjugate pairs of points in $\Qcone(\F_{q^2})$ and each pair imposes two linear conditions on curves in $\Cone \to V^{\mathrm{con}} \cong \A^{15}$. Thus
\[
S_0(\Cone_{[2]}) = (q^4 - q) q^{13}.
\]  
Similarly,
\[
S_0(\Cone_{[2,1]}) = (q^4 -q)(q^2+q) q^{12} \mbox{ \ \ \ and \ \ \ } S_0(\Cone_{[3]}) = (q^6 + q^3 - q^2 - q) q^{12}.
\]
Being singular at a specified point imposes $3$ linear conditions unless $\lambda = [2,1]$ and the singularity is at the $\F_q$-rational marked point. The singularity needs to lie outside every line containing at least two of the marked points, because otherwise the fiber must contain that line and hence the vertex of the cone. Thus
\[
S_1(\Cone_{[2]}) = (q^4 - q^2)s_1(X) q^{10} \mbox{ \ \ and \ \ } S_1(\Cone_{[3]}) = (q^6 - q^2) s_1(X) q^9,
\]
while
\[
S_1(\Cone_{[2,1]}) = (q^4 - q^2) (q^2+q)s_1(\mathrm{pt})q^{10} +  (q^4-q-(q^2-q))(q^2+q-1) s_1(X) q^9.
\]
The computations of $S_d(\Cone_\lambda)$ for $d = 2$ and $3$ are likewise closely analogous to the computations for $\lambda = [1^2]$ and $[1^3]$ in the proof of Proposition~\ref{prop:ConeSieve}, using the Hasse--Weil sieve and considering cases according to how many of the singularities coincide with marked points and how many are coplanar with all three singularities.
\end{proof}

\begin{proposition} \label{prop:NonsplitEquiv}
The sieve terms $S_d(\Cns_\lambda)$  for $\lambda \in \{ [2], [2,1], [3] \}$ and $d \leq 3$ up to $o(q^6)$ are 
\[
\begin{array}{|l||l|}
\hline
S_0(\Cns_{[2]}) &  (q^4+q^2)\p_{13}
\Tstrut \\ [2.5 pt] \hline
S_1(\Cns_{[2]}) &  -(q^2+1)(q^4+q^2) \p_{10}
\Tstrut \\ [2.5 pt] \hline
S_2(\Cns_{[2]}) &  -q^{13} -q^{12} -q^{11} - q^{10}
\Tstrut \\ [2.5 pt] \hline
S_3(\Cns_{[2]}) &  4q^{12} + 3q^{11} + 5q^{10} + 5q^9 + 3q^8 + 3q^7 + 3q^6 + \cdots 
\Tstrut \\ [2.5 pt] \hline
\end{array}
\]
\[
\begin{array}{|l||l|}
\hline
S_0(\Cns_{[2,1]}) &  (q^2+1)(q^4+q^2) \p_{12}
\Tstrut \\ [2.5 pt] \hline
S_1(\Cns_{[2,1]}) &  -(q^2+1)(q^4+q^2)\p_{10} -(q^2+1)(q^6+q^4) \p_9
\Tstrut \\ [2.5 pt] \hline
S_2(\Cns_{[2,1]}) &  q^{15} - q^{14} +q^{13} -2q^{12} -q^{11} -q^{10} - q^9
\Tstrut \\ [2.5 pt] \hline
S_3(\Cns_{[2,1]}) &  4q^{13} + 6q^{12} + 8q^{11} +10q^{10} +7q^9 + 7q^8 +4q^7 + 4q^6 + \cdots 
\Tstrut \\ [2.5 pt] \hline
\end{array}
\]
\[
\begin{array}{|l||l|}
\hline
S_0(\Cns_{[3]}) &  q^{18} + q^{17} + q^{16} + q^{15} + \cdots 
\Tstrut \\ [2.5 pt] \hline
S_1(\Cns_{[3]}) &  -q^{17} -q^{16} -2q^{15} -2q^{14} -q^{13} -q^{12} + q^7 + q^6 + \cdots 
\Tstrut \\ [2.5 pt] \hline
S_2(\Cns_{[3]}) & 0
\Tstrut \\ [2.5 pt] \hline
S_3(\Cns_{[3]}) &  q^{13} + 2q^{11} + 2q^{10} + q^8 -2q^7 -2q^6 + \cdots
\Tstrut \\ [2.5 pt] \hline
\end{array}
 \ .
\]
\end{proposition}

\begin{proof}
The proof is similar to that of Proposition~\ref{prop:NonsplitSieve} for $n = 2$ and $3$.  The surface $\Qns$ has $q^4 + q^2$ ordered conjugate pairs of points over $\F_{q^2}$ and $q^6 - q^2$ ordered conjugate triples over $\F_{q^3}$.  All such collections of marked points impose independent linear conditions on curves in the linear series $\Cns \to \P^{15}$.  Thus
\[
S_0(\Cns_{[2]}) = (q^4 + q^2) \p_{13}, \quad S_0(\Cns_{[2,1]}) = (q^4 + q^2)(q^2+1) \p_{12}, \quad \mbox{ and } \quad S_0(\Cns_{[3]}) = (q^6 -q^2) \p_{12},
\]
where $\p_m = \# \P^m(\F_q)$. To compute $S_1$, note that a singularity imposes $3$ additional linear conditions unless it coincides with a marked point, in which case it imposes only $2$ conditions.

For $S_2$, the number of linear conditions imposed by a pair of singularities depends only on the number of singularities that coincide with marked points.  The vanishing of $S_2(\Cns_{[3]})$ is immediate; none of the singularities can coincide with marked points and so this sieve coefficient is divisible by $s_2(\Qns) = 0$.

For $S_3$, the three singularities span a unique hyperplane section. We consider cases according to whether this hyperplane section is smooth or singular, how many of the marked points lie in this hyperplane, and how many marked points coincide with singularities.   
\end{proof}

\begin{proposition} \label{prop:SplitEquiv}
The sieve terms $S_d(\Csplit_\lambda)$ for $\lambda \in \{ [2], [2,1], [3] \}$ and $d \leq 3$ up to $o(q^6)$ are 
\[
\begin{array}{|l||l|}
\hline
S_0(\Csplit_{[2]}) &  q^{17} + q^{16} +2q^{15} + \cdots 
\Tstrut \\ [2.5 pt] \hline
S_1(\Csplit_{[2]}) & -q^{16} -3q^{15} -5q^{14} -5q^{13} -2q^{12} + \cdots 
\Tstrut \\ [2.5 pt] \hline
S_2(\Csplit_{[2]}) &  2 q^{14}+5 q^{13}+11 q^{12}+5 q^{11}+q^{10}-2 q^6 + \dots
\Tstrut \\ [2.5 pt] \hline
S_3(\Csplit_{[2]}) &  -2 q^{12}-7 q^{11}-7 q^{10}-q^9+q^8+q^7+3 q^6 +\dots
\Tstrut \\ [2.5 pt] \hline
\end{array}
\]
\[
\begin{array}{|l||l|}
\hline
S_0(\Csplit_{[2,1]}) &  q^{18} + 3q^{17} + 5q^{16} + 5q^{15} + 2q^{14} + \cdots 
\Tstrut \\ [2.5 pt] \hline
S_1(\Csplit_{[2,1]}) &  -q^{17}-6 q^{16}-16 q^{15}-22 q^{14}-15 q^{13}-4 q^{12} + q^7 + 5 q^6 + \dots
\Tstrut \\ [2.5 pt] \hline
S_2(\Csplit_{[2,1]}) &  3 q^{15}+17 q^{14}+35 q^{13}+32 q^{12}+11 q^{11}-q^{10}-q^9
-2 q^8-8 q^7-18 q^6 + \dots
\Tstrut \\ [2.5 pt] \hline
S_3(\Csplit_{[2,1]}) &   -8 q^{13}-24 q^{12}-26 q^{11}-10 q^{10}+3 q^9+7 q^8+12 q^7+18 q^6 +\dots
\Tstrut \\ [2.5 pt] \hline
\end{array}
\]
\[
\begin{array}{|l||l|}
\hline
S_0(\Csplit_{[3]}) &  q^{18} +q^{17} + q^{16} + 3q^{15} + 2q^{14} + \cdots 
\Tstrut \\ [2.5 pt] \hline
S_1(\Csplit_{[3]}) &  -q^{17} - 3 q^{16} - 6 q^{15} - 10 q^{14} - 9 q^{13} - 3 q^{12} + q^7 + 3 q^6 + \cdots 
\Tstrut \\ [2.5 pt] \hline
S_2(\Csplit_{[3]}) & 2 q^{15}+8 q^{14}+16 q^{13}+16 q^{12}+6 q^{11}-2 q^8
-4 q^7-6 q^6 + \dots
\Tstrut \\ [2.5 pt]  \hline
S_3(\Csplit_{[3]}) &  -3 q^{13}-12 q^{12}-14 q^{11}-4 q^{10}+2 q^9+3 q^8+4 q^7+6 q^6 +\dots
\Tstrut \\ [2.5 pt] \hline
\end{array}
 \ .
\]
\end{proposition}

\begin{proof}
The proof is similar to those of Proposition~\ref{prop:SplitSieve} for $n =2$ and $3$, and Proposition~\ref{prop:NonsplitEquiv}. There are $q^4 + q^2 - 2q$ ordered conjugate pairs of points in $\Qsplit(\F_{q^2})$ and $q^6 + 2q^3 -q^2 - 2q$ ordered conjugate triples in $\Qsplit(\F_{q^3})$. Each such configuration of marked points imposes independent linear conditions on $\Csplit \to \P^{15}$.  Thus
\[
\resizebox{.99\hsize}{!}{
$
S_0(\Csplit_{[2]}) = (q^4 + q^2 - 2q) \p_{13}, \quad S_0(\Csplit_{[2,1]}) = (q^4 + q^2 - 2q)(q^2+2q+1) \p_{12}, \quad S_0(\Csplit_{[3]}) = (q^6 + 2q^3 -q^2 - 2q) \p_{12}.
$}
\]
A singularity imposes 3 additional linear conditions unless it coincides with a marked point or is collinear with 3 marked points. The exceptional cases, where the singularity imposes only 2 new linear conditions, cannot occur for $\lambda = [2]$, and hence
\[
S_1(\Csplit_{[2]}) = s_1(\Qsplit) (q^4 + q^2 - 2q) \p_{10}.
\]
For $\lambda = {[3]}$, the singularity cannot coincide with a marked point, but it could be collinear with all 3 marked points.  There are $2(q+1)$ lines, each of which contains $q^3-q$ ordered conjugate triples over $\F_{q^3}$. Therefore
\[
S_1(\Csplit_{[3]}) = s_1(\Qsplit)(q^6+2q^3-q^2-2q) \p_9 + 2(q+1) s_1(\P^1) (q^3-q) (\p_{10}-\p_9).
\]
The computation of $S_1(\Csplit_{[2,1]})$ is similar, accounting in addition for the cases where the singularity coincides with the $\F_q$-rational marked point.

The remaining computations of $S_2$ and $S_3$ are similar to the corresponding computations in the proof of Proposition~\ref{prop:SplitSieve}; substantial bookkeeping is required to keep track of all of the cases, but each case presents no new difficulties.
\end{proof}

\subsection{Proof of Theorem~\ref{thm:equiv-point-counts}}
This proof is analogous to the proof of Theorems~\ref{thm:polynomial} and \ref{thm:openpolynomial} in \S\ref{sec:mainproof}.  
We begin with $n=2$. Putting together the information from Propositions~\ref{prop:ConeEquiv}, \ref{prop:NonsplitEquiv} and \ref{prop:SplitEquiv} we get the approximate $\SS_2$-equivariant counts, 
\[
\begin{array}{|l||l|}
\hline
\mbox{\raisebox{-2pt}{$\#^{\SS_2} \Cone_2(\F_q)/\# \mathrm{Aut}(\Qcone)$}}  &\mbox{\raisebox{-2pt}{$(q^{10}+q^9-q^8-q^7+o(q^{\frac{11}{2}}) ) s_2+$}}\\ [2 pt] &
+ (q^9-q^8-q^7+q^6+o(q^{\frac{11}{2}}) )  s_{1^2}
\Tstrut \\ \hline
\mbox{\raisebox{-2pt}{$\#^{\SS_2} \Cns_2(\F_q)/\# \mathrm{Aut}(\Qns)$}}   & \mbox{\raisebox{-2pt}{$\frac{1}{2} (q^{11}-q^7 +o(q^{\frac{11}{2}}) ) s_2+$}}  \\ [2 pt] &
+\frac{1}{2} (-q^9+q^8-q^6+o(q^{\frac{11}{2}}) )  s_{1^2}
\Tstrut \\ [2 pt] \hline
\mbox{\raisebox{-1.5pt}{$\#^{\SS_2}\Csplit_2(\F_q)/\# \mathrm{Aut}(\Qsplit)$}}  & \mbox{\raisebox{-2pt}{$\frac{1}{2} (q^{11}+2q^{10}+2q^9-4q^8-q^7+o(q^{\frac{11}{2}}) ) s_2+$}} \\ [2 pt] &
+  \frac{1}{2} (2q^{10}+q^9-q^8-2q^7-q^6+o(q^{\frac{11}{2}}) )  s_{1^2}
\Tstrut \\ [1.5 pt] \hline
\end{array}
\]

\noindent Adding the $\SS_2$-equivariant count of the hyperelliptic locus in Section~\ref{sec:equivarianthypers}, we get 
$\#^{\SS_2}\cM_{4,2}(\F_q)$ up to $o(q^{\frac{11}{2}})$. We add $\#^{\SS_2} \partial \cM_{4,2}(\F_q)$, see below, computed using the Getzler-Kapranov formula as in the proof of Proposition~\ref{prop:boundary-count} to get $\#^{\SS_2} \ocM_{4,2}(\F_q)$ up to $o(q^{\frac{11}{2}})$. Then $\#^{\SS_2} \ocM_{4,2}(\F_q)$ is determined by Poincar\'e duality, as in Remark~\ref{rem:PD}, and $\#^{\SS_2} \cM_{4,2}(\F_q)$ is obtained by subtracting the contribution of the boundary. 

For $n=3$ analogous computations show:
\[
\begin{array}{|l||l|}
\hline
\mbox{\raisebox{-1pt}{$\#^{\SS_3} \Cone_3(\F_q)/\# \mathrm{Aut}(\Qcone)$}}  &
\mbox{\raisebox{-1pt}{$(q^{11}+q^{10}-q^9-q^8+q^7+*\,q^6+o(q^6) )s_3+$}} \\ [1 pt] &
+(q^{10}-q^9-3q^8+2q^7+*\,q^6+o(q^6) )s_{2,1}+\\&
+(-q^8+q^7+*\,q^6+o(q^6) )s_{1^3}
\Tstrut \\ \hline
\mbox{\raisebox{-1pt}{$\#^{\SS_3} \Cns_3(\F_q)/\# \mathrm{Aut}(\Qns)$}}   & 
\mbox{\raisebox{-1pt}{$\frac{1}{2}(q^{12}-q^8+*\,q^6+o(q^6) )s_3+$}} \\ [1 pt] &
+\frac{1}{2}(-q^{10}+q^9+q^8-2q^7+*\,q^6+o(q^6) )s_{2,1}+\\&
+\frac{1}{2}(-q^{10}+q^8-2q^7+*\,q^6+o(q^6) )s_{1^3}
\Tstrut \\ \hline
\mbox{\raisebox{-1pt}{$\#^{\SS_3}\Csplit_3(\F_q)/\# \mathrm{Aut}(\Qsplit)$}}  & 
\mbox{\raisebox{-1pt}{$\frac{1}{2}(q^{12}+2q^{11}+2q^{10}-4q^9-3q^8+2q^7+*\,q^6+o(q^6) )s_3+$}} \\ [1 pt] &
+\frac{1}{2}(2q^{11}+3q^{10}-3q^9-3q^8+*\,q^6+o(q^6) )s_{2,1}+\\&
+\frac{1}{2}(q^{10}-q^8+*\,q^6+o(q^6) )s_{1^3}
\Tstrut \\ [1 pt] \hline
\end{array}
\]
\noindent Here, each $*$ denotes an undetermined integer. Together with $\#^{\SS_3}\cH_{4,3}(\F_q)$, these computations determine $\#^{\SS_3} \cM_{4,3}(\F_q)$ up to $o(q^6)$. Adding $\#^{\SS_3} \partial \cM_{4,3}(\F_q)$, then using the polynomiality of $\#^{\SS_3} \ocM_{4,3}(\F_q)$ together with Poincar\'e duality, and subtracting $\#^{\SS_3} \partial \cM_{4,3}(\F_q)$ we get 
\begin{multline*}
\#^{\SS_3} \cM_{4,3}(\F_q)=
(q^{12}+2q^{11}+3q^{10}-2q^9-4q^8+2q^7+C_1 \, q^6-q^3+2q^2)s_3+\\+(q^{11}+2q^{10}-q^9-4q
^8+q^7+C_2 \, q^6-4q^3-1)s_{2,1}+(-q^8+C_3\, q^6-2q^3+2q+1)s_{1^3}
\end{multline*}
for some as yet unknown integers $C_1,C_2$ and $C_3$.
Evaluating at $q=1$ gives the $\SS_3$-equivariant Euler characteristic $\chi^{\SS_3}( \cM_{4,3})=2s_3-6s_{2,1}$, see \cite{Gorsky14}. This shows that $C_1=-1,C_2=0$ and $C_3=0$, and completes the proof. 

We end by recording the equivariant point counts of the boundaries, computed using the Getzler-Kapranov formula.
\[
\begin{array}{|l||l|}
\hline
\mbox{\raisebox{-1pt}{$\#^{\SS_2}(\partial \cM_{4,2})(\F_q)$}} &  \mbox{\raisebox{-1pt}{$(7q^{10}+52q^9+222q^8+563q^7+901q^6+901q^5+561q^4+221q^3+56q^2+$}}\\&\mbox{\raisebox{-.75pt}{$+9q+1)s_2+
(q^{10}+20q^9+99q^8+279q^7+461q^6+461q^5+277q^4+$}}\\
&\mbox{\raisebox{-.5pt}{$+100q^3+22q^2+2q)s_{1^2}$}}
\Tstrut \\ \hline
\mbox{\raisebox{-1pt}{$\#^{\SS_3}(\partial \cM_{4,3})(\F_q)$}} &  \mbox{\raisebox{-1 pt}{$(9q^{11}+84q^{10}+426q^9+1351q^8+2692q^7+3415q^6+2694q^5+1347q^4+$}}\\
&\mbox{\raisebox{-.75pt}{$+425q^3+85q^2+11q+1)s_3+(4q^{11}+56q^{10}+350q^9+1224q^8+2577q^7+$}}\\
&\mbox{\raisebox{-.5pt}{$+3304q^6+2578q^5+1220q^4+353q^3+58q^2+5q+1)s_{2,1}+(4q^{10}+46q^9+$}}\\
&\mbox{\raisebox{-.25pt}{$+191q^8+446q^7+583q^6
+446q^5+190q^4+48q^3+4q^2-2q-1)s_{1^3}$}}
\Tstrut \\ \hline
\end{array}
\]

\subsection{Counts in terms of local systems} \label{sec:localsystems}
The universal curve $\pi \colon \cM_{4,1} \to \cM_4$ gives rise to the $\ell$-adic local system $\V:=R^1 \pi_* \Q_{\ell}$ whose stalk over $[C] \in \cM_4$ is $H^1(C,\Q_{\ell})$. For every $\lambda= \lambda_1, \lambda_2, \lambda_3, \lambda_4$, with $\lambda_1\geq\lambda_2\geq\lambda_3 \geq\lambda_4\geq 0$, we get, by applying Schur functors, an induced local system $\mathbb{V}_{\lambda}$ from the irreducible representation of $\mathrm{GSp}(8)$ with the corresponding highest weight. The trace of Frobenius on the compactly supported Euler characteristic of $\cM_4$ with coefficients in $\V_{\lambda}$ for all $\lambda$ such that $|\lambda|=\lambda_1+\ldots+\lambda_4 \leq n$ gives equivalent information to computing 
$\#^{\SS_m} \cM_{4,m}$ for all $0 \leq m \leq n$, see \cite{Getzler99} and \cite{Bergstrom08}. For simplicity of notation, we omit trailing zeros and write, for instance, $\V_1 := \V_{1,0,0,0}$. Using this we can show the following.
\begin{theorem}  \label{thm:M4loc}
For any $q$ we have, 
\[
\begin{array}{|l||l|}
\hline
\mathrm{Tr}(F_q,
    H^\bullet_c(\cM_{4},\V_{1}))&  q^7+q^2
    \Tstrut \\ \hline
\mathrm{Tr}(F_q,
    H^\bullet_c(\cM_{4},\V_{2}))& 0 
        \Tstrut \\ \hline
\mathrm{Tr}(F_q,
    H^\bullet_c(\cM_{4},\V_{1,1})) &  q^9-q^8-q^7-q^2
\Tstrut \\ \hline
\mathrm{Tr}(F_q,
    H^\bullet_c(\cM_{4},\V_{3})) &  q^3-2q-1
\Tstrut \\ \hline
\mathrm{Tr}(F_q,
    H^\bullet_c(\cM_{4},\V_{2,1})) &  q^8-q^7-q^4+2q^3+1
\Tstrut \\ \hline
\mathrm{Tr}(F_q,
    H^\bullet_c(\cM_{4},\V_{1,1,1})) &  -q^{10}-q^9+q^8+q^6-q^4-q^2
 \Tstrut \\ \hline
\end{array}
\]
\end{theorem}
\noindent  
Note that the equalities of Theorem~\ref{thm:M4loc} can be translated into equalities of Euler characteristics with values in the Grothendieck group of either $\ell$-adic Galois representations or Hodge structures  as in Remark~\ref{rem:eulergalois}.

\bibliographystyle{amsalpha}
\bibliography{Genus4Counting}  

\providecommand{\bysame}{\leavevmode\hbox to3em{\hrulefill}\thinspace}
\providecommand{\MR}{\relax\ifhmode\unskip\space\fi MR }
\providecommand{\MRhref}[2]{%
  \href{http://www.ams.org/mathscinet-getitem?mr=#1}{#2}
}
\providecommand{\href}[2]{#2}
\begin{thebibliography}{BCGY21}

\bibitem[AC98]{ArbarelloCornalba98}
E.~Arbarello and M.~Cornalba, \emph{Calculating cohomology groups of moduli
  spaces of curves via algebraic geometry}, Inst. Hautes \'{E}tudes Sci. Publ.
  Math. (1998), no.~88, 97--127.

\bibitem[AV04]{ArsieVistoli04}
A.~Arsie and A.~Vistoli, \emph{Stacks of cyclic covers of projective spaces},
  Compos. Math. \textbf{140} (2004), no.~3, 647--666.

\bibitem[BCGY21]{BibbyChanGadishYun21}
C.~Bibby, M.~Chan, N.~Gadish, and C.~Yun, \emph{Homology representations of
  compactified configuration spaces on graphs applied to $\mathcal{M}_{2,n}$},
  arXiv:2109.03302, to appear in Exp. Math., 2021.

\bibitem[Beh91]{BehrendThesis}
K.~Behrend, \emph{The {L}efschetz trace formula for the moduli stack of
  principal bundles}, ProQuest LLC, Ann Arbor, MI, 1991, Thesis
  (Ph.D.)--University of California, Berkeley.

\bibitem[Beh93]{Behrend93}
\bysame, \emph{The {L}efschetz trace formula for algebraic stacks}, Invent.
  Math. \textbf{112} (1993), no.~1, 127--149.

\bibitem[Ber08]{Bergstrom08}
J.~Bergstr\"{o}m, \emph{Cohomology of moduli spaces of curves of genus three
  via point counts}, J. Reine Angew. Math. \textbf{622} (2008), 155--187.

\bibitem[Ber09]{Bergstrom09}
\bysame, \emph{Equivariant counts of points of the moduli spaces of pointed
  hyperelliptic curves}, Doc. Math. \textbf{14} (2009), 259--296.

\bibitem[BG01]{BrockGranville01}
B.~Brock and A.~Granville, \emph{More points than expected on curves over
  finite field extensions}, Finite Fields Appl. \textbf{7} (2001), no.~1,
  70--91, Dedicated to Professor Chao Ko on the occasion of his 90th birthday.

\bibitem[BP00]{BoggiPikaart00}
M.~Boggi and M.~Pikaart, \emph{Galois covers of moduli of curves}, Compositio
  Math. \textbf{120} (2000), no.~2, 171--191.

\bibitem[BT07]{BergstromTommasi07}
J.~Bergstr{\"o}m and O.~Tommasi, \emph{The rational cohomology of
  {$\overline{\mathcal M}_4$}}, Math. Ann. \textbf{338} (2007), no.~1,
  207--239.

\bibitem[CFP12]{ChurchFarbPutman12}
T.~Church, B.~Farb, and A.~Putman, \emph{The rational cohomology of the mapping
  class group vanishes in its virtual cohomological dimension}, Int. Math. Res.
  Not. IMRN (2012), no.~21, 5025--5030.

\bibitem[Cha21]{Chan21}
M.~Chan, \emph{Topology of the tropical moduli spaces ${M}_{2,n}$}, Beitr.
  Algebra Geom. (2021), 1--25.

\bibitem[CL19]{ChenevierLannes19}
G.~Chenevier and J.~Lannes, \emph{Automorphic forms and even unimodular
  lattices}, Ergebnisse der Mathematik und ihrer Grenzgebiete. 3. Folge. A
  Series of Modern Surveys in Mathematics, vol.~69, Springer, 2019.

\bibitem[CL22]{CanningLarson22}
S.~Canning and H.~Larson, \emph{On the {C}how and cohomology rings of moduli
  spaces of stable curves}, preprint arXiv:2208.02357v1, 2022.

\bibitem[CR15]{ChenevierRenard15}
G.~Chenevier and D.~Renard, \emph{Level one algebraic cusp forms of classical
  groups of small rank}, Mem. Amer. Math. Soc. \textbf{237} (2015), no.~1121,
  v+122.

\bibitem[DSvZ21]{admcycles}
V.~Delecroix, J.~Schmitt, and J.~van Zelm, \emph{admcycles---a {S}age package
  for calculations in the tautological ring of the moduli space of stable
  curves}, J. Softw. Algebra Geom. \textbf{11} (2021), no.~1, 89--112.

\bibitem[FHR21]{FlorenceHoffmannReichstein21}
M.~Florence, N.~Hoffmann, and Z.~Reichstein, \emph{On the rationality problem
  for forms of moduli spaces of stable marked curves of positive genus}, Ann.
  Sc. Norm. Super. Pisa Cl. Sci. (5) \textbf{22} (2021), no.~3, 1091--1104.

\bibitem[Get98]{Getzler98}
E.~Getzler, \emph{The semi-classical approximation for modular operads}, Comm.
  Math. Phys. \textbf{194} (1998), no.~2, 481--492.

\bibitem[Get99]{Getzler99}
\bysame, \emph{Resolving mixed {H}odge modules on configuration spaces}, Duke
  Math. J. \textbf{96} (1999), no.~1, 175--203.

\bibitem[GK98]{GetzlerKapranov98}
E.~Getzler and M.~Kapranov, \emph{Modular operads}, Compositio Math.
  \textbf{110} (1998), no.~1, 65--126.

\bibitem[Gor14]{Gorsky14}
E.~Gorsky, \emph{The equivariant {E}uler characteristic of moduli spaces of
  curves}, Adv. Math. \textbf{250} (2014), 588--595.

\bibitem[Kee92]{Keel92}
S.~Keel, \emph{Intersection theory of moduli space of stable {$n$}-pointed
  curves of genus zero}, Trans. Amer. Math. Soc. \textbf{330} (1992), no.~2,
  545--574.

\bibitem[KL02]{KisinLehrer02}
M.~Kisin and G.~I. Lehrer, \emph{Equivariant {P}oincar\'{e} polynomials and
  counting points over finite fields}, J. Algebra \textbf{247} (2002), no.~2,
  435--451.

\bibitem[Log03]{Logan03}
A.~Logan, \emph{The {K}odaira dimension of moduli spaces of curves with marked
  points}, Amer. J. Math. \textbf{125} (2003), no.~1, 105--138.

\bibitem[Mes86]{Mestre86}
J.-F. Mestre, \emph{Formules explicites et minorations de conducteurs de
  vari\'{e}t\'{e}s alg\'{e}briques}, Compositio Math. \textbf{58} (1986),
  no.~2, 209--232.

\bibitem[MSS13]{MoritaSakasaiSuzuki13}
S.~Morita, T.~Sakasai, and M.~Suzuki, \emph{Abelianizations of derivation {L}ie
  algebras of the free associative algebra and the free {L}ie algebra}, Duke
  Math. J. \textbf{162} (2013), no.~5, 965--1002.

\bibitem[Pet17]{Petersen17}
D.~Petersen, \emph{A spectral sequence for stratified spaces and configuration
  spaces of points}, Geom. Topol. \textbf{21} (2017), no.~4, 2527--2555.

\bibitem[Pik95]{Pikaart95}
M.~Pikaart, \emph{An orbifold partition of {$\overline M{}^n_g$}}, The moduli
  space of curves ({T}exel {I}sland, 1994), Progr. Math., vol. 129,
  Birkh\"{a}user Boston, Boston, MA, 1995, pp.~467--482.

\bibitem[PW21]{PayneWillwacher21}
S.~Payne and T.~Willwacher, \emph{The weight 2 compactly supported cohomology
  of moduli spaces of curves}, preprint arXiv:2110.05711v1, 2021.

\bibitem[Ste]{SF}
J.~Stembridge, \textsf{SF}, a {M}aple package for symmetric functions. This
  package is available at
  \texttt{http://www.math.lsa.umich.edu/\~{}jrs/maple.html}.

\bibitem[Tom05]{Tommasi05}
O.~Tommasi, \emph{Rational cohomology of the moduli space of genus 4 curves},
  Compos. Math. \textbf{141} (2005), no.~2, 359--384.

\bibitem[vdBE05]{vandenBogaartEdixhoven05}
T.~van~den Bogaart and B.~Edixhoven, \emph{Algebraic stacks whose number of
  points over finite fields is a polynomial}, Number fields and function
  fields---two parallel worlds, Progr. Math., vol. 239, Birkh\"{a}user Boston,
  Boston, MA, 2005, pp.~39--49.

\bibitem[vR16]{vanRooij}
S.~B. van Rooij, \emph{Het tellen van krommen op het product van twee
  projectieve lijnen over een eindig lichaam}, Thesis (B.Sc.)--Utrecht
  University, 2016.

\bibitem[VW15]{VakilWood15}
R.~Vakil and M.~Wood, \emph{Discriminants in the {G}rothendieck ring}, Duke
  Math. J. \textbf{164} (2015), no.~6, 1139--1185.

\bibitem[Wen20]{Wennink20}
T.~Wennink, \emph{Counting the number of trigonal curves of genus 5 over finite
  fields}, Geom. Dedicata \textbf{208} (2020), 31--48.

\bibitem[Won22]{Wong22}
Y.~M. Wong, \emph{An algorithm to compute fundamental classes of spin
  components of strata of differentials}, preprint arXiv:22111.16061v1, to
  appear in Int. Math. Res. Not. (IMRN), 2022.

\bibitem[Xar20]{Xarles20}
X.~Xarles, \emph{A census of all genus 4 curves over the field with 2
  elements}, preprint arXiv:2007.07822, 2020.

\end{thebibliography}
\end{document}